\documentclass[11pt,titlepage]{article}
\usepackage[margin=1.2in]{geometry}
\usepackage{amsmath, amsbsy, amsthm, array, mathrsfs}
\usepackage{graphics}
\usepackage{physics}
\usepackage{epsfig, amssymb,latexsym,verbatim}
\usepackage{graphicx}
\usepackage{color}
\usepackage{dsfont, bbm}
\usepackage{relsize}
\usepackage{subcaption}

\newcommand{\R}{\mathbb{R}}
\newcommand{\mbs}{\mathbb{S}}

\newcommand{\noin}{\noindent}
\newcommand{\bee}{\begin{eqnarray*}}
\newcommand{\ene}{\end{eqnarray*}}
\newcommand{\bec}{\begin{center}}
\newcommand{\enc}{\end{center}}
\newcommand{\be}{\begin{equation}}
\newcommand{\ee}{\end{equation}}

\newcommand{\ep}{\varepsilon}
\newcommand{\mb}{\mathbf}
\newcommand{\bs}{\boldsymbol}
\newcommand{\tb}{\textbf}
\newcommand{\pend}{$\blacksquare$}
\newcommand{\vs}{\vskip 3mm}
\newcommand{\bi}{\begin{itemize}}
\newcommand{\ei}{\end{itemize}}
\begin{document}

\title{\LARGE  Robustness of  the deepest projection regression functional
 \\[4ex]
}

\author{ {\sc Yijun Zuo}\\[2ex]
         {\small {\em  Department of Statistics and Probability,  Michigan State University} }\\
         {\small East Lansing, MI 48824, USA} \\
         {\small zuo@msu.edu}\\[6ex]
     }
 \date{\today}
\maketitle

\vskip 3mm
{\small

\begin{abstract}

Depth notions in regression have been systematically proposed and examined in Zuo (2018).  One of the prominent advantages of the notion of depth is that it can be directly utilized to introduce median-type deepest estimating functionals (or estimators in the case of empirical distributions) for location or regression parameters in a multi-dimensional setting. \vskip 3mm

Regression depth shares the advantage. Depth induced deepest estimating functionals are expected to inherit desirable and inherent robustness properties ( e.g. bounded maximum bias and influence function and high breakdown point) as their univariate location counterpart does. Investigating and verifying the robustness of the deepest projection estimating functional (in terms of maximum bias, asymptotic and finite sample breakdown point, and influence function)  is the major goal of this article.
\vskip 3mm
 It turns out that the deepest projection estimating functional possesses a bounded influence function and the best possible asymptotic breakdown point as well as the best finite sample breakdown point with robust choice of its univariate regression and scale component.
\bigskip
\par

\noindent{\bf Key words and phrase:} Depth,  linear regression, deepest regression estimating functionals,
maximum bias, breakdown point, influence function, robustness.
\bigskip
\par
\noindent {\bf Running title:} Robustness of deepest regression functional.
\end{abstract}
}

\section{Introduction}

Consider a general linear regression model
\begin{eqnarray}
y&=&\mathbf{x}'\boldsymbol{\beta}+{{e}}, \label{eqn.model}
\end{eqnarray}
where $y$ and ${e}$ are univariate random variables, $'$ denotes the transpose of a vector, and random vector $\mathbf{x}=(x_1,\cdots, x_p)'$ and unknown  parameter $\boldsymbol{\beta}$ 
are in $\R^p$, the error ${e}$ has distribution $F_{{e}}$ and 
the random vector $\mathbf{x}$ has distribution
$F_{\mathbf{x}}$. Note that this general model includes the special case with an intercept term. For example,
if $\bs{\beta}=(\beta_1, \bs{\beta_2}')'$ and $x_1=1$, then one has $y=\beta_1+\mb{x_2}'\bs{\beta_2}+{e}$, where $\mb{x_2}=(x_2,\cdots, x_p) \in \R^{p-1}$.
If one denotes $\mb{w}=(1,\mb{x_2}')'$, then $y=\mb{w}'\bs{\beta}+{e}$. We use this model or (\ref{eqn.model}) interchangably depending on the context.
Denote by $F_{(y,~\mathbf{x})}$ the joint distribution of  $y$ and $\mathbf{x}$ under the model (\ref{eqn.model}). 
\vs
Let $T(\cdot )$ be a $\R^p$-valued  estimating functional for $\boldsymbol{\beta}$, defined on the set ${\mathcal G}$ of distributions on $\R^{p+1}$.  $T$ is called \emph{Fisher consistent} for  $\boldsymbol{\beta}$  if $T( F_{(y,~\mathbf{x})})=\boldsymbol{\beta}_0$ for the true parameter $\bs{\beta}_0\in \R^p$ of the model and for $F_{(y,~\mathbf{x})} \in {\mathcal G_1}\subset {\mathcal G}$,
each member of ${\mathcal G_1}$ possesses some common attributes.
Additional desirable properties of a regression functional $T(\cdot)$ are \emph{regression}, \emph{scale}, and \emph{affine} equivariant. That is,
$ T(F_{(y+\mathbf{x}'b, ~\mathbf{x})})=T(F_{(y,~\mathbf{x})})+b, \forall ~b \in \R^p;$~~
$ T(F_{(sy, ~\mathbf{x})})=s T(F_{(y, ~\mathbf{x})}), \forall ~s \in \R;$ 
$ T(F_{(y,~ A'\mathbf{x})})=A^{-1}T(F_{(y,~ \mathbf{x})}),\forall \mbox{~nonsingular  $A\in \R^{p\times p}$} $;  respectively.
Namely, $T(\cdot)$ does not depend on the underlying coordinate system and measurement scale.
\vskip 3mm
The classical regression estimating functional
 is the least square (LS)
 functional. It meets all the desired properties above and is ``optimal" if $F_e$ is normal (Huber (1972)). But it is extremely sensitive to a slight deviation from the normality assumption.
 Alternatives include the least absolute deviation functional, and quantile regression (Koenker and Bassett (1978)) were posed.  But in terms of asymptotic breakdown point (ABP) robustness, they are no better than the traditional LS functional (all have $0\%$ ABP). Estimating functionals with higher ABP  were consequently proposed. Among them, the least median squares estimator (Rousseeuw (1984)) is the most famous one.
It has the highest ABP ($50\%$) but suffers a slow convergence rate (cubic root) (Davies (1989 and Kim and Pollard (1990)) and a instability drawback (Figure 3.2 of Seber and Lee (2003)).
\vs
Robust estimating functionals with high ABP and root $n$ convergence rate were subsequently advanced. Among many of them is the regression depth (RD) induced deepest regression estimating functional (Rousseeuw and Hubert (1999) (RH99)) ($T^*_{RD}$). The latter has an ABP $1/3$ (Van Aelst and Rousseeuw (2000) (VAR00)) and root $n$ consistency (Bai and He (1999)).
\vs
One of the prominent advantages of depth notion
is that it can be directly employed to introduce median-type deepest estimating functionals (or estimators in the empirical case) for the location or regression parameter in a multi-dimensional setting based on a general min-max stratagem.
The most outstanding feature of the univariate median is its exceptional robustness. Indeed, it has the best possible finite sample breakdown point (FSBP) (among all location equivariant estimators, see Donoho (1982)) and the minimum maximum bias (MB) (if the underlying distribution has a unimodal symmetric density, see Huber (1964)).
\vs
The functional in RH99 ($T^*_{RD}$) holds desired properties, its ABP ($1/3$) is lower than the highest ($1/2$) though.
The deepest projection estimating functional ($T^*_{PRD}$) induced from projection regression depth (PRD) in Zuo (2018)(Z18) overcomes this. It has the best ABP with a root n consistency ((Z18)) as well. $T^*_{PRD}$ is closely related to the bias-robust estimates (P-estimates) of Marrona and Yohai (1993) (MY93). In fact, it is a modified version of the latter, achieving the scale-invariance (see Section 2).
\vs

MY93 investigated the robustness of P-estimates, provided an upper bound of their MB, but their influence function (IF) and FSBP had not been explicitly established in the last quarter of century.
Establishing a  MB upper bound for $T^*_{PRD}$ and discovering its IF
and revealing its exact FSBP are three main objectives of this article.
\vskip 3mm

The rest of the article is organized as follows.
 Section 2 introduces the $T^*_{PRD}$.
 Section 3 is devoted to the establishment of MB, IF and FSBP of $T^*_{PRD}$. Section 4 addresses the computation issues of the deepest regression estimators, and presents data examples to illustrate the performance (in terms of robustness) of the regression lines of the LS, the $T^*_{RD}$ and the $T^*_{PRD}$, and carries out some simulations to investigate the finite-sample relative efficiency of $T^*_{RD}$ and the $T^*_{PRD}$.
Brief concluding remarks  end the article in Section 5.
\vskip 3mm

\section{Maximum projection regression depth functionals} \label{T^*-sec}

Let us first recall  the projection regression depth and its induced deepest estimating functionals defined in Z18.
\vskip 3mm
Assume that $T$ is a {univariate} regression estimating functional
 which satisfies\\[1ex]
\noindent
(\tb{A1}) {regression}, {scale} and {affine equivariant}, that is, \vskip 2mm

$T(F_{(y+x{b},~{x})})=T(F_{(y,~{x})})+{b}$, ~$\forall ~{b}\in \R$, and \vskip 2mm

$T(F_{(sy,~{x})})=sT(F_{(y,~{x})})$, ~$\forall ~s  \in\R$, and\vskip 2mm

 $T(F_{(y,~a{x})})=a^{-1}T(F_{(y,~{x})})$,~ $\forall ~a (\neq 0)\in\R$.\vskip 2mm

respectively, where $x, y\in \R$ are random variables (r.v.s). Throughout the lower case $x$ is in $R$ while bold $\mb{x}$ is a vector.\vskip 3mm
\vs
\noindent
 Let $S$ be a positive scale estimating functional such that \vskip 3mm

 \noindent
(\tb{A2}) $S(F_{sz+b})=|s|S({F_z})$ \text{for any r.v. $z \in \R$ and scalar $b, s \in \R$}, that is, $S$ is \emph{scale equivariant} and \emph{location invariant}.
 \vskip 3mm

Equipped with a pair of $T$ and $S$, we can introduce a corresponding projection based multiple regression estimating functional.
  Define
\be\mbox{UF}_\mb{v}(\boldsymbol{\beta}; ~F_{(y,~\mathbf{x})}, T):= |T(F_{(y-\mathbf{x}'\boldsymbol{\beta},~\mathbf{x}'\mb{v})})|/S(F_{y}), \label{eqn.uf}
\ee
which represents
 unfitness of $\boldsymbol{\beta}$ at $F_{(y,~\mathbf{x})}$  w.r.t. $T$ along the $\mb{v}\in \mbs^{p-1}:=\{u: \|u\|=1, u\in\R^p\}$.
 If $T$ is a \emph{Fisher consistent} regression estimating functional, then $T(F_{(y-\mathbf{x}'\boldsymbol{\beta}_0,~\mathbf{x}'\mb{v})})=0$
 for some $\bs{\beta}_0$ (the true parameter of the model) and $\forall ~\mb{v}\in \mbs^{p-1}$.
 Then, overall one expects $|T|$ to be small and close to zero for a candidate $\bs{\beta}$, independent of the choice of $\mb{v}$ and $\mathbf{x}'\mb{v}$. {The magnitude of $|T|$ measures the unfitness of $\bs{\beta}$ along the $\mb{v}$}.
 Taking the supremum over all $\mb{v}\in\mbs^{p-1}$,  yields
 \be
 \mbox{UF}(\boldsymbol{\beta};~F_{(y,~\mathbf{x})},T)= \sup_{\|\mb{v}\|=1}\mbox{UF}_\mb{v}(\boldsymbol{\beta}; ~F_{(y,~\mathbf{x})}, T),\label{eqn.UF}\ee
 the\emph{ unfitness} of $\boldsymbol{\beta}$ at $F_{(y,~\mathbf{x})}$  w.r.t. 
 $T$. 
 Now applying the min-max scheme, we obtain the \emph{projection regression estimating functional} (also denoted by $T^*_{\text{PRD}}$) w.r.t. the pair $(T,S)$
 \begin{eqnarray}
 T^*(F_{(y,~\mathbf{x})},T)&=&\arg\!\min_{\boldsymbol{\beta}\in \R^p}\mbox{UF}(\boldsymbol{\beta}; ~F_{(y,~\mathbf{x})}, T)  \label{eqn.T*}\\[1ex]
 &=&\arg\!\max_{\bs{\beta}\in\R^p}\text{PRD}\left(\boldsymbol{\beta}; ~F_{(y,~\mathbf{x})}, T\right),\nonumber
 \end{eqnarray} 
 where, the \emph{projection regression depth} (PRD) function is defined as
 \be
 \text{PRD}\left(\boldsymbol{\beta}; ~F_{(y,~\mathbf{x})}, T\right)=\left(1+\mbox{UF}\big(\boldsymbol{\beta}; ~F_{(y,~\mathbf{x})}, T\big)\right)^{-1},
 \label{eqn.PRD}
 \ee
\vs
\noindent

\noindent
 \textbf{Remarks 2.1}\vskip 3mm
(I)~$\mbox{UF}(\boldsymbol{\beta};~F_{(y,~\mathbf{x})}, T)$ or $\mbox{UF}(\boldsymbol{\beta};~F_{(y,~\mathbf{x})})$ corresponds to outlyingness $O(x, F)$, and  $T^*(F_{(y,~\mathbf{x})})$ corresponds
to the projection median functional $PM(F)$ in location setting (see Zuo (2003)). Note that in (\ref{eqn.uf}), (\ref{eqn.UF}) and (\ref{eqn.T*}), we have suppressed the scale $S$ since it does not involve $\mb{v}$ and is nominal. 
Sometimes we also suppress $T$ for convenience.\vskip 3mm

 A similar $T^*$  was first introduced and studied in MY93, where it was called P1-estimate (denote it by $T_{P1}$, see (\ref{eqn.MY93})). However, they are  different.
 The definition of $T^*$ here is \emph{different from} $T_{P1}$ of MY93, the latter multiplies by $S(F_{\mb{v}'\mathbf{x}})$ instead of dividing by $S(F_y)$ in $\mbox{UF}_v(\boldsymbol{\beta}; ~F_{(y,~\mathbf{x})}, T)$ here.  Furthermore, MY93 did not talk about the ``unfitness" (or ``depth").
 Corresponding to (\ref{eqn.uf}) here,  they instead defined the following
 $$
 A(\boldsymbol{\beta}, \mb{v}) =| T(F_{(y- \boldsymbol{\beta}'\mathbf{x}, ~\mb{v}'\mathbf{x})})| S(F_{\mb{v}'\mathbf{x}}),
 $$
 where $\mb{v}, \boldsymbol{\beta} \in \R^p$. Their P1-estimate is defined as
 \be
 T_{P1}=\arg \min_{\boldsymbol{\beta}\in \R^p}\sup_{\|\mb{v}\|=1}A(\boldsymbol{\beta}, \mb{v}).\label{eqn.MY93}
 \ee

(II) It is readily seen that $A(\bs{\beta}, \mb{v})$ is \emph{not scale invariant} whereas $\mbox{UF}_{\mb{v}}(\boldsymbol{\beta};~F_{(y,~\mathbf{x})}, T)$ is. $T^*$ is regression, scale, and affine equivariant. 
 \vs
 (III)  Examples of $T$ include mean, quantile, and median( Med),  and location functionals in Wu and Zuo (2009) (WZ09). Examples of $S$ include standard deviation functional, the median absolute deviations functional (MAD), and scale functionals in WZ08. Hereafter we write $\text{Med}(Z)$ rather than
  $\text{Med}(F_Z)$.
 For the special choice of $T$  and $S$ in (\ref{eqn.uf}) such as
 \bee
 T(F_{(y-\mb{x}'\bs{\beta},~ \mb{x'}\mb{v})})&=&\text{Med}_{\mb{x'}\mb{v}\neq 0}\big(\frac{y-\mb{x}'\bs{\beta}}{\mb{x'}\mb{v}}\big),\\[1ex]\label{spesific-T.eqn}
 S(F_y)&=& \text{MAD}(F_y), \label{S.eqn}
 \ene
we have
 \be
 \text{UF}(\bs{\beta}; F_{(y,~\mb{x})})=\sup_{\|\mb{v}\|=1}\Big|\text{Med}_{\mb{x'}\mb{v}\neq 0}\big(\frac{y-\mb{x}'\bs{\beta}}{\mb{x'}\mb{v}}\big)\Big|\bigg/ \text{MAD}(F_y),
 \ee
 and
 \be
 \text{PRD}\left(\bs{\beta}; ~F_{(y,~\mathbf{x})}\right)=\inf_{\|\mb{v}\|=1,\mb{x'}\mb{v}\neq 0}
 \frac{\text{MAD}(F_y)}{\text{MAD}(F_y)+\Big|\text{Med}\big(\frac{y-\mb{x}'\bs{\beta}}{\mb{x'}\mb{v}}\big)\Big|}. \label{special-PRD.eqn}
 \ee\vskip 3mm

A special case of PRD above (the empirical case) is closely related to the so-called ``centrality" in Hubert, Rousseeuw, and Van Aelst (2001) (HRVA01). In the definition of the latter, nevertheless, all the term of ``MAD$(\cdot)$" on the RHS  of (\ref{special-PRD.eqn}) is divided by $\text{Med}|\mb{x'}\mb{v}|$.
 \hfill \pend
\section{Robustness of the deepest projection regression functional}

One of the main purposes of seeking the maximum depth estimating functional in regression is for the robustness consideration, since the classical LS functional is notorious sensitive to the deviation from the model assumptions (normality assumption) and to the contamination. On the other hand, a maximum depth estimating functional could be regarded as a median-type functional in regression. The latter in location is well-known for its exceptional robustness.
 Do the maximum projection depth estimating functionals inherit the inherent robustness properties of the location counterpart (and w.r.t. what types of robustness measure)? \vskip 3mm
\subsection{Maximum bias}
 For a given distribution $F \in\R^d$ (hereafter $F\in\R^d$ really means that $F$ is defined on $\R^d$) and an $\ep > 0$, the version of $F$ contaminated by an $\ep$ amount of an \emph{arbitrary distribution} $G\in\R^d$ is
denoted by $F(\ep, G) = (1 -\ep)F + \ep G$ (an $\ep$ amount deviation from the assumed $F$). Here it is assumed that $\ep \leq 1/2$, otherwise $F(\ep, G)=G((1-\ep), F)$, and one can't distinguish which one is contaminated by which one. The \emph{maximum bias} of a given
general functional L under an $\ep$ amount of contamination at $F$ is defined as (see Hampel,
Ronchetti, Rousseeuw and Stahel (1986) (HRRS86))
\bee \mbox{MB}(\ep; L, F) &=& \sup_{G \in \R^d} \|L (F (\ep, G)) - L (F )\|, \ene
where $\mbox{MB}(\ep; L,F)$ is the maximum deviation (bias) of $L$ under an $\ep$ amount
of contamination at $F$ and it mainly measures the global robustness of $L$.
For a given $L$ at $F$, it is desirable that $\mbox{MB}(\ep; L, F)$ is bounded for an $\ep(\leq 1/2)$ as large as possible. 
\vskip 3mm

The minimum amount $\ep^*$ of contamination at $F$ which leads to an unbounded
$\mbox{MB}(\ep; L,F)$ is called the \emph{asymptotic breakdown point} (ABP) of $L$ at $F$,
$\ep^* (L, F) = \inf\{\ep :\mbox{MB}(\ep; L,F) = \infty\}$.
\vskip 3mm

\noindent
For a given $F=F_{(y,~\mb{x})}\in\R^{p+1}$, write $F_{(\mb{v}, \bs{\beta})}:=F_{(y-\mb{x}'\mb{\bs\beta},~\mb{x}'\mb{v})}$ for $\mb{v} \in \mbs^{p-1}$ and a given $\bs{\beta}\in\R^p$.
Let $F_y$ be the marginal distribution based on $y\in \R$.
For
the univariate regression (and scale) estimating
functional $T$ (and $S$) in Section \ref{T^*-sec} and an $\ep>0$, define
\bee
B_T(\ep; T, F) &=& \inf_{\bs{\beta}\in\R^{p}} \sup_{G \in \R^2, \|\mb{v}\|=1}|T(F_{(\mb{v},\bs{\beta})}(\ep, G))|,\\[1ex]
C(\ep; T,F)&=& \sup_{G \in \R^2, \|\mb{v}\|=1}|T(F_{(\mb{v},\mb{0})})(\ep,G))|,\\[1ex]
B_S(\ep; S, F)&=& \sup_{G\in \R}|S(F_y(\ep, G))|,~~~~~~ 
b(\ep; S, F)~=~ \inf_{G\in \R}|S(F_y(\ep, G))|.
\ene
\vskip 3mm

\noindent
\tb{Proposition 3.1} For a given pair $(T,S)$,~ 
 $F=F_{(y,~\mb{x})}$, and an $\ep>0$,  assume that $T(F_{(\mb{v},\bs{0})})=0$,  $\forall~\mb{v}\in \mbs^{(p-1)}$, and $b(\ep;S, F)>0$ and $B_S(\ep; S, F)<\infty$. Then for
$T^*( F_{(y,~\mb{x})}, T)$  in (\ref{eqn.T*})
$$\mbox{MB}(\ep; T^*, F)\leq ~  B_T(\ep; T, F)+C(\ep; T, F).
$$

\vskip 3mm

\noindent
\textbf{Proof:} By regression equivariance of the $T^*$ (see (II) of Remarks 2.1), 
assume (w.l.o.g) that 
$T^*(F)=0$. Then
$$\mbox{MB}(\ep; T^*, F)=\sup_{G\in\R^{p+1}}\|T^*(F(\ep, G)\|.
$$
For the given $F$ and a given $G\in\R^{p+1}$, denote $\bs{\beta}^*(\ep, G):=T^*(F(\ep, G))$ and $F(\ep,G)=F_{\mb{z}^*}$ and $\mb{z}^*=(y^*, {\mb{x}^*}')'\in\R^{p+1}$. Then we need to show that
$$\sup_{G\in\R^{p+1}}\|\bs{\beta}^*(\ep, G)\| \leq  B_T(\ep; T, F)+C (\ep; T, F).$$
For the given $G\in\R^{p+1}$ and  $F$, by (\ref{eqn.uf}), (\ref{eqn.UF}), and (\ref{eqn.T*}),  we have
\bee
\bs{\beta}^*(\ep, G)&=&\arg\! \min_{\bs{\beta}\in\R^{p}} \mbox{UF}(F(\ep,G); \bs{\beta},T)\\[1ex]
&=&\arg\! \min_{\bs{\beta}\in\R^{p}} \! \sup_{\|v\|=1} \mbox{UF}_{\mb{v}}(F(\ep,G); \bs{\beta},T)
\ene
Assume that  $\bs{\beta}^*(\ep, G)\neq 0$. 
Write $\bs{\beta}^*$ for $\bs{\beta}^*(\ep, G)$ and let $v^*=\bs{\beta}^*/\|\bs{\beta}^*\|$,  then we have by \tb{(A1)} given in Section \ref{T^*-sec}
$$|T(F_{(y^*-(\mb{x}^*)'\bs{\beta^*},~(\mb{x}^*)'\mb{v^*})})|=|T(F_{(y^*,~(\mb{x}^*)'\mb{v^*})})-\|\bs{\beta}^*\||.$$
If $\|\bs{\beta}^*\|\leq \sup_{\|\mb{v}\|=1}|T(F_{(y^*,~(\mb{x}^*)'\mb{v})}|$ for every given $G \in R^{p+1}$,  
then $\|\bs{\beta}^*\|\leq C(\ep; T, F) $, we already have the desired result. Otherwise, we have
 for any given $\bs{\beta} \in\R^{p}$
\bee\sup_{\|\mb{v}\|=1}|T(F_{(y^*-(\mb{x}^*)'\bs{\beta},~(\mb{x}^*)'\mb{v})})|&\geq&
|T(F_{(y^*-(\mb{x}^*)'\bs{\beta^*},~(\mb{x}^*)'\mb{v^*})})|\\[1ex]
&\geq& \|\bs{\beta}^*\|-|T(F_{(y^*,~(\mb{x}^*)'\mb{v^*})}|\\[1ex]
&\geq& \|\bs{\beta}^*\|-\sup_{\|\mb{v}\|=1}|T(F_{(y^*,~(\mb{x}^*)'\mb{v})}|,
 \ene
Therefore, we have for the given $G\in\R^{p+1}$ and $F$ and $\ep$ and the given $\bs{\beta}\in\R^p$
\bee
\|\bs{\beta}^*(\ep, G)\|&\leq& \sup_{\|\mb{v}\|=1}|T(F_{(y^*-(\mb{x}^*)'\bs{\beta},~(\mb{x}^*)'\mb{v})})|+\sup_{\|\mb{v}\|=1}|T(F_{(y^*,~(\mb{x}^*)'\mb{v})})|\\[1ex]
&\leq& \sup_{G \in\R^2,\|\mb{v}\|=1}|T(F_{(\mb{v},\bs{\beta})}(\ep,G))|+\sup_{G \in\R^2, \|\mb{v}\|=1}|T(F_{(\mb{v},\bs{0})}(\ep,G))|. 
\ene
Taking the infimum over $\bs{\beta} \in\R^p$ and then supremum over $G\in \R^{(p+1)}$ in both sides immediately yields the desired result.
This completes the proof.\hfill \pend
\vskip 3mm

\noindent
\tb{Remarks 3.1}
\vs
{(I)} The assumption \tb{(A0)}: $T(F_{(\mb{v},0)})=T(F_{(y,~\mb{x}'\mb{v})})=0$ for $\mb{v}\in\mbs^{p-1}$ is equivalent to the Fisher-consistency of $T$ or $F_{(y,~\mb{x})}$ is T-symmetric about $\mb{0} \in \R^p$.
 $F_{(y,~\mb{x})}$   is \emph{T-symmetric} about a $\bs{\beta}_0$ iff 

\be \tb{(C0)}:~~~ T\big(F_{(y-\mb{x'}\bs{\beta}_0,~\mb{x}'\mb{v})}\big)=0,~\forall ~\mb{v}\in\mbs^{p-1},
 \label{T-symm.eqn}\ee
and it holds for a wide range of distributions $F_{(y,~\mb{x})}$ and $T$.
For example, if the univariate
functional $T$ is the mean functional, then this becomes the classical assumption in regression when $\bs{\beta}_0$ is the true parameter of the model:
the conditional expectation of the error term ${e}$ (which is assumed to be independent of $\mb{x}$)  given $\mb{x}$  is zero, i.e. $$ \tb{(C1)}:~~~~ E(F_{y-\mb{x'}\bs{\beta_0}} \big|_{\mb{x}=\mb{x_0}})=E(F_{(y-\mb{x'}\bs{\beta_0},\mb{x'}\mb{v}} )=0,~\forall ~\mb{x_0}\in\R^{p-1}, \mb{v}\in \mbs^{p-1}.\hspace*{8mm} $$
\vs
\noin
\tb{(A0)}, however, is not indispensable in the proof but for the neatness of the upper bound and of the expression for $B_T(\ep; T,F)$. Adding $\sup_{\|\mb{v}\|=1}|T(F_{(y,\mb{x}'\mb{v})})|$ to the RHS of the upper bound and using the regular \emph{deviation} definition for $B_T(\ep;T,F)$, the proposition holds without \tb{(A0)}.
 \vskip 3mm
{(II)} An upper bound for their P-estimates was also given in MY93 (Theorem 3.3).
The two upper bounds are quite different due to the definition of $T^*$ is different from P-estimates.
\vskip 3mm

{(III)} The conditions on $S$ in the proposition are typically satisfied by common scale functionals such as
MAD or scale functionals in WZ08.
The term $C(\ep; T,F)$ in the Proposition is typically bounded for $T$ (such as quantile functionals
or functionals in WZ09).
\vskip 3mm
{(IV)} The maximum projection regression depth functional $T^*$ has a bounded maximum bias as long as that is true for the $T$, and $S$ does not breakdown (for a scale functional, its ABP is defined as
$\ep^*(S, F)=\min\{\ep: B_S(\ep; S,F)+b(\ep; S, F)^{-1}=\infty \}$). Furthermore, the MB upper bound of $T^*$  depends entirely on that of the $T$ as long as $S$ does not breakdown. The Proposition also reveals the ABP of
$T^*$ as summarized in the following. ~\hfill \pend
\vskip 3mm
\noindent
\tb{Corollary 3.1} Under the same assumptions of Proposition 3.1, 
we have
\begin{itemize}
\item[](i)
$\ep^*(T^*, F)\geq\min\big \{\ep^*(T, F); \ep^{*}(S, F)  \big\}$.\vskip 3mm
\item[]\hspace*{-10mm} if $(T, S)=(\text{Med},~\text{MAD})$ 
 then\vskip 2mm
\item[](ii) $\ep^*(T^*, F) = 1/2$
\end{itemize}
\vskip 3mm

\noindent
\textbf{Proof:} \vskip 3mm
(i) is trivial.\vskip 3mm

(ii) follows from the standard ABP results of Med and MAD (see e.g. HRRS86) and the upper bound of ABP for any regression equivariant functional (see Theorem 3.1 of Davies (1993) and of Davies and Gather (2005)).  \hfill \pend
\vskip 3mm

\noin
\tb{Remarks 3.2}\vskip 3mm
(I) If the choice for $T$ and $S$ is $(\text{Med},~\text{MAD})$, then $T^*$ can have an ABP as high as $1/2$. HRVA01 reported  their most central regression estimator $T^c_r$ (in Theorem 8) has a $50\%$ breakdown point without any rigorous treatment. $T^c_r$, however, is  slightly different from $T^*(F_n)$ here, see Remarks 2.1.

\vskip 3mm
(II) The ABP of the deepest regression functional of RH 99
 has been inventively studied in VAR00 and is $1/3$, 
while the ABP of the classical LS functional is $0$.
\hfill\pend
\vskip 3mm
When $(T, S)$ is $(\text{Med},~\text{MAD})$ , then the general bounds involved in Proposition 3.1 could be concretized and
specified as shown in the following. Furthermore, one also could construct a lower bound for the maximum bias of
$T^*$ in (\ref{eqn.T*}).
\vskip 3mm
 First we need some notations. Write $q(\ep)=1/(2(1-\ep))$ for a given $0<\ep <1/2$. Denote $m_i(Z, c,\ep)$ for quantiles such that  $m_1(Z,c,\ep)=F_{|Z-c|}^{-1}(1-q(\ep))$, $m_2(Z, c,\ep)=F_{|Z-c|}^{-1}(q(\ep))$ for a random variable $Z\in\R$ any scalar $c\in \R$ . \vskip 3mm
\noin
\tb{Proposition 3.2}~
Let $T(F_{(y-\bs{x'}\bs{\beta},~\mb{x}'\mb{v})})=\text{Med}\big(\frac{y-\bs{x'}\bs{\beta}}{\mb{x}'\mb{v}}\big)$ ($\mb{x}'\mb{v}\neq 0$ a.s.),~ $S(F_y)= \text{MAD}(F_y)$.
 Assume that
$1^o$) $F_{(y, \bs{x})}$ is $T$-symmetric about a $\bs{\beta_0}$ which is the true parameter of model (\ref{eqn.model});
$2^o$) $F_{{e}}$ has a symmetric, decreasing in $|x|$ density $f(x)$; 
$3^o$)  $F_{\mb{x}'\mb{v}}$ is the same $\forall~\mb{v}\in\mbs^{p-1}$;
$4^o$) $e$ and $\mb{x}$ are independent.
Then,  for the $T^*$ in (\ref{eqn.T*}), the given $F=F_{(y,~\mb{x})}$, any  $0<\ep<1/2$, \vskip 3mm
\begin{itemize}
\item[] (i) $T^*$ is Fisher-consistent. That is, $T^*(F, T)=\bs{\beta}_0$, under $1^o)$;\vspace*{-0mm}
\item[] (ii) $B_S(\ep; S,F)=c$,~$b(\ep; S,F)=d$, under $1^o)$--$2^o)$;  $B_T(\ep; T, F)=b$, under  $3^o$);
\item[] (iii) $\displaystyle b \leq \mbox{MB}(\ep; T^*, F) \leq  b+C(\ep;T, F)=2b$, under $1^o)$--$4^o)$;
\end{itemize}
where 
 $b=J^{-1}\big(q(\ep)\big)$,
 $c=m_2(y,a1,\ep)$, $d=m_1(y, b1,\ep)$, $a1=F_{|y|}^{-1}(1-q(\ep))$, $b1=F_{|y|}^{-1}(q(\ep))$. 
 All quantiles is assumed to exist uniquely,
 $J$ is the distribution of $y/\mb{x}'\mb{v}$,$\mb{v}\in \mbs^{p-1}$. 
\vskip 3mm
\noindent
To prove the statements above, we need the following result given in Zuo, Cui, and Young (2004) (ZCY04).
\vskip 3mm
\noin
\tb{Lemma 3.1}~ Suppose that $A = F^{-1}(1-q(\ep))$ and $B = F^{-1}( q(\ep))$ exist uniquely
for $X\in\R$ with $F:=F_X$ and $0 < \ep < 1/2$.  Let $\delta_x$ denote the point-mass probability measure at $x \in \R$. Then for any distribution $G \in \R$ and point $x$,
\begin{itemize}
\item[] (L-i) $A \leq {\text{Med}}(F(\ep,G)) \leq B$, 
 ~~~(L-ii) ${\text{Med}}(F(\ep, \delta_x )) = {\text{Med}}\{A, B, x\}$,\vskip 3mm
\item[] (L-iii) $m1\big(X,{\text{Med}}(F(\ep, G)), \ep\big) \leq {\text{MAD}}(F(\ep, G)) \leq m2\big(X,{\text{Med}}(F(\ep, G)), \ep\big)$, \vskip 3mm
\item[] (L-iv) ${\text{MAD}}(F(\ep, \delta_x)) = {\text{Med}}\bigg\{m1\big(X,\text{Med}(F (\ep, \delta_x)), \ep\big), ~ |x - {\text{Med}}(F (\ep, \delta_x ))|,\\
    \hspace*{48mm} ~ m2\big(X,{\text{Med}}(F (\ep, \delta_x)), \ep\big)\bigg\}$.
\end{itemize}
where $\text{Med}$ is applied to distributions as well as discrete points. \hfill  \pend

\vskip 3mm
\noindent
\tb{Proof of Proposition 3.2} \vskip 3mm
(i)  The given condition (assumption) guarantees that $T$ is Fisher-consistent at $F_{(y,\mb{x})}$, that is, for any $\mb{v}\in\mbs^{p-1}$ $$T\big( F_{(y-\mb{x}'\bs{\beta_0},~\mb{x}'\mb{v})}\big)=0.$$
Both (\ref{eqn.uf}) and (\ref{eqn.UF}) are equal to zero. That is, $ \text{UF}(\bs{\beta}_0; F_{(y,\mb{x})}, T)=0.$ Therefore, $\bs{\beta}_0$ attains the minimum possible value of $\text{UF}(\bs{\beta}; F_{(y,\mb{x})}, T)$ for any $\bs{\beta}\in \R^{p}$, which further means that $T^*(F_{(y,~\mb{x})},T)=\bs{\beta}_0$.
By the equivalence of $T^*$ (see Remarks 2.1), assume w.l.o.g. that $\bs{\beta}_0=\mb{0}$.
\vskip 3mm

(ii) We need the maximum bias bounds on Med and MAD. Some of them have been already established in Lemma A.2 of ZCY04 (cited above in Lemma 3.1). 

\vskip 3mm

Note that when $\bs{\beta_0}=0$, $y$ has the same distribution as ${e}$, $m_i(y,c,\ep)$ is nonincreasing in $c$ for $c>0$, 
the bounds for $S$ follow directly from this fact, coupled with  (L-iii) and (L-i). 
\vskip 3mm
We have to establish the bound for $T$.  Note that
\be B_T(\ep; T, F)=\inf_{\bs{\beta}\in\R^p}\sup_{G\in R^2, \|\mb{v}\|=1}|T(F_{(\mb{v},\bs{\beta})}(\ep; G))|.
\label{B-eqn}
\ee
To invoke (L-i) of the Lemma 3.1, we need to first determine the $B$ in (L-i) for the distribution of $Z:=(y-\mb{x}'\bs{\beta})/(\mb{x}'\mb{v})$ for a given $\bs{\beta}\in\R^p$ and a $\mb{v}\in\mbs^{p-1}$. Note that
\bee
Z=\frac{y-\mb{x}'\bs{\beta}}{\mb{x}'\mb{v}}&=&
\frac{y-\mb{x}'(\bs{\beta}-(\bs{\beta}'\mb{v})\mb{v})-\mb{x}'\mb{v}(\bs{\beta}'\mb{v})}{\mb{x}'\mb{v}}\\[1ex]
&:=&\frac{y-\mb{x}'\bs{\alpha}_{(\bs{\beta},\mb{v})}} {\mb{x}'\mb{v}} -\bs{\beta}'\mb{v}\\[1ex]
&:=& Z1-\bs{\beta}'\mb{v}.
\ene
For convenience we suppress the dependency of $Z$ and $Z1$ on $\bs{\beta}$ and $\mb{v}$. Note that $\|\bs{\alpha}_{(\bs{\beta},\mb{v})}\|= \|\bs{\beta}-(\bs{\beta}'\mb{v})\mb{v}\|=(\|\bs{\beta}\|^2-(\bs{\beta}'\mb{v})^2)^{1/2}$ and $\bs{\alpha}'_{(\bs{\beta},\mb{v})}\mb{v}=0$. 
It is readily seen that
$F_Z(z)=F_{Z1}(z+\bs{\beta}'\mb{v})$ and hence that $F_Z^{-1}(p)=F_{Z1}^{-1}(p)-\bs{\beta}'\mb{v}$ for any $p\in (0,1)$. \vskip 3mm
Now denote the distribution of $(y-\mb{x}'\bs{\alpha})/\mb{x}'\mb{v}$ with
$\|\bs{\alpha}\|=r$ and $\bs{\alpha}'\mb{v}=0$ by $J_{r}$ for any  $\mb{v}\in\mbs^{p-1}$.
Hence $F_{Z1}=J_r$ with $r=(\|\bs{\beta}\|^2-(\bs{\beta}'\mb{v})^2)^{1/2}$ for any $\mb{v}\in\mbs^{p-1}$ and a given $\bs{\beta}\in\R^p$. 
\vskip 3mm
In the light of Lemma 3.1,  \be B_T(\ep;T, F)=\inf_{\bs{\beta}\in\R^p}\sup_{\|v\|=1}\left| F^{-1}_{Z1}(q(\ep))-\bs{\beta}'\mb{v}\right|. \label{B.eqn}\ee
On the other hand,
\bee\inf_{\bs{\beta}\in\R^p}\sup_{\|v\|=1}\left| F^{-1}_{Z1}(q(\ep))-\bs{\beta}'\mb{v}\right|
&\geq& \inf_{\bs{\beta}\in\R^p}\left| J^{-1}_{0}(q(\ep))-\|\bs{\beta}\|\right|\\[.5ex]
&=&\inf_{\bs{\beta}\in\R^p} \left|\|\bs{\beta}\|- J^{-1}_{0}(q(\ep))\right|
= J^{-1}_{0}(q(\ep)),
\ene
where the first inequality follows from the consideration of a special $\mb{v}=\bs{\beta}/\|\bs{\beta}\|$ for $\bs{\beta}\neq 0$ and $3^o$), the second equality is due to the fact that $J^{-1}_{0}(q(\ep))$ has nothing to do with $\bs{\beta}$. Therefore, by picking $\bs{\beta}=\mb{0}$ on the RHS of (\ref{B.eqn}), its LHS attains its lower bound.
That is, $B_T(\ep; T,F)=J_0^{-1}(q(\ep))$ which is the same as $b$ since when $r=0$,  $J_r$ is the same as $J$ distributionally. 
\vskip 3mm

(iii) In virtue of (ii) above, one part of the RHS inequality has already been established in Proposition 3.1. But we still need to show that $C(\ep;T,F)=b$. This, however, follows in a straightforward manner from the definition of $C(\ep; T,F)$ and the proof in (ii) above
(with $\bs{\beta}=\mb{0}$ in this case).
\vskip 3mm
 We need to show the LHS lower bound for $\mbox{MB}(\ep; T^*, F)$. We adapt the idea of Huber (1981) (page 74-75). Note that for a given $\mb{v}\in\mbs^{p-1}$ by $4^o$) $$F_{\mb{v}}(y,{z}):=F_{(y,~\mb{x}'\mb{v})}(y,{z})=F_{(y-\bs{x}'\bs{\beta}_0,~\mb{x}'\mb{v})}(y,{z})
=F_{{e}}(y)F_{\mb{x}'\mb{v}}({z}), ~\text{for~} y, ~z\in \R.$$

Assume that $\bs{x}\neq \bs{0}$, otherwise, our discussion reduces to Huber (1981) (page 74-75), our conclusion holds true. Assume, w.l.o.g., that the first component of $\bs{x}$, $x_1\neq 0$.
Construct two functions:
$$
F^+_{\mb{v}}(y, {z})= (1-\ep)\big[ F_{{e}}(y)\mb{I}_{y\leq ax_1}+ F_{{e}}(y-2ax_1)\mb{I}_{y>ax_1}\big]F_{\mb{x}'\mb{v}}({z}),
$$
$$
F^{-}_{\mb{v}}(y,{z})=F^{+}_{(y,\mb{x})}(y+2ax_1,{z}),
$$
where $a=J^{-1}(q(\ep))$ is the $q(\ep)$th quantile of $y/\bs{x'}\bs{v_0} (=y/x_1)$ with $\bs{v_0}=(1,0,\cdots,0)'\in\R^p$. It is now not difficult to verify that the two functions above are distribution functions over $\R^{2}$ and belong to $F_{\mb{v}}(\ep; G)$ for some $G \in \R^{2}$(because both keep $(1-\ep)$ part of $F_{\mb{v}}(y, {z})$).\vskip 3mm
Assume that for some the random vector $(y^*, \mb{x})$, $F_{(y^*, ~\mb{x}'\mb{v})}=F_{\mb{v}}^{+}$. (note that vector $\mb{x}$ is unchanged due to the construction).
Then one has $F_{\mb{v}}^{-}=F_{(y^*+2ax_1,~\mb{x}'\mb{v})}=F_{(y^*+\mb{x}'\bs{\eta},~\mb{x}'\mb{v})}$ with $\bs{\eta}=(2a,0,\cdots,0)'\in\R^p$
\vskip 3mm
Denote the first coordinate of $T^*(F)$ as $T^*_1 (F)$. Then  by the equivariance of $T^*$, we see that
$T_1^*(F_{\mb{v_0}}^{+})-T^*_1(F_{\mb{v_0}}^{-})=-2 a,$
which implies 
\bee
2a&\leq &\sup_{\|\mb{v}\|=1}|T_1^*(F_{\mb{v}}^{+})-T_1^*(F_{\mb{v}}^{-})|\\[1ex]
&\leq& 2\sup_{G\in\R^{2}, \|\mb{v}\|=1}\|T^*(F_{\mb{v}}(\ep; G))\|. 
\ene
\vskip 3mm
\noindent
Note that $a=b$. This completes the entire proof.
\hfill \pend
\vskip 3mm

\noindent
\tb{Remarks 3.3}\vskip 3mm

(I) Part (i) of the Proposition holds as long as $T$ is $T$-symmetric about a $\bs{\beta}_0 \in \R^{p}$. That is, $T$ is not necessarily to be the Med functional. Furthermore, $S$ plays no role in the verification process, that is, any scale estimating functional will work.
 Likewise,  the lower bound in (iii) holds true for any $T$ 
 and $S$. 
The (Med, MAD) choice is just the classical one.\vskip 3mm

(II) The assumption that  $F_{{e}}$ has a symmetric density $f(x)$ which is decreasing in $|x|$ is common and  typically required in the literature (see, e.g., MY93, Theorem 3.5). It guarantees that the construction of the two functions are indeed
distribution functions in the proof of (iii) (actually it guarantees that the probability mass covered by both $F_{{e}}(y)\mb{I}_{y\leq ax_1}$ and $F_{{e}}(y-2ax_1)\mb{I}_{y>ax_1}$ are $q(\ep)$, therefore guarantees the success of the construction).
\vskip 3mm
(III) The assumption $3^o$), that is, $F_{\mb{x}'\mb{v}}$ is the same for any $v\in\mbs^{p-1}$
holds if (i) $(y, \mb{x}'\mb{v})$ is spherically distributed about the origin or (ii) if $\mb{x}$ is spherically distributed about the origin. 
(ii) was assumed in Theorem 3.5 of MY93. However, in the light of the equivalence of $T^*$, the spherical symmetry could be relaxed to elliptical symmetry. 
\vskip 3mm

(IV) In many cases, the maximum bias is attained by a point-mass distribution, that is, $\mbox{MB}(\ep; T^*,F) =
\sup_{x\in\R^d}\| T^* (F (\ep, \delta_x )) - T^* (F )\|$ (see Huber (1964), Martin, Yohai and Zamar
(1989), Chen and Tyler (2002) and Adrover and Yohai (2002)). 
The upper bound in (iii) also appeared in MY93 ((a) of Theorem 4.1). Where it was shown attainable by  a variant of their P1-estimate (different from $T^*$ here) under the point-mass contamination $\delta_x$ and when
  $X$ is spherical distributed.
\hfill \pend
\vskip 3mm

Maximum bias and ABP are global robustness measure and depict the global robust perspectives of the underlying functional. Now we will focus on the local robustness  of $T^*$ via its influence function.\vskip 3mm
\subsection{Influence function}
The \emph{influence function} (IF) of a functional $T$ at a given point $x \in \R^d$ for a given $F$ is defined
as
$$
\text{IF}(x; T,F) = \lim_ {\ep\to 0^+} \frac{T (F (\ep, \delta_x)) - T (F )}{\ep},$$
where $\delta_x$ is the point-mass probability measure at $x \in \R^d$ , and the \emph{gross error
sensitivity} of $T$ at $F$ is then defined as (in HRRS86)
$$
\gamma^* (T , F ) = \sup_{x\in\R^d}\|
\text{IF}(x; T,F)\|. $$

The function $\text{IF}(x; T,F)$ describes the relative effect (influence) on $T$ of an
infinitesimal point-mass contamination at $x$ and measures the local robustness
of $T$ . The function $\gamma^*(T , F )$ is the maximum relative effect on $T$ of an
infinitesimal point-mass contamination and measures the global as well as local
robustness of $T$.  It is desirable that a regression estimating functional has a bounded influence function
and especially a bounded gross-error sensitivity. This, however, does not hold for an arbitrary regression estimating functional, especially for the  classical least squares functional. Now we investigate this for $T^*$ in (\ref{eqn.T*}). \vskip 3mm

 For the sake of simplicity,  we will assume  below that
$\mb{x}$ is spherically distributed, i.e. the distribution of $\mb{x'}\mb{v}$ is the same for any $\mb{v}\in\mbs^{p-1}$. The result and the discussion, however, can be trivially extended to cover the case that $\mb{x}$ is elliptically distributed, in the light of  the equivalence of $T^*$ (see Remarks 2.1) and Proposition 1 of VAR00.\vs
Denote $\mb{z}:=(y,\mb{x})$, $F(y,\mb{s}):=F_{\mb{z}}(y, \mb{s})$. Consider the point-mass $\ep$ contamination of $F_{(y,\mb{x})}$  at $\delta_{\mb{z}}$:  $F_{(y,\mb{x})}(\ep; \delta_{\mb{z}})=(1-\ep)F_{(y,\mb{x})}+\ep\delta_{(y_0, \mb{x_0})}$, where $\mb{x_0}=(x_{01},x_{02},\cdots,x_{0p})'\in\R^p$ and $\mb{x_0}\neq0$.
Denote $z_0:=y_0/x_{01}$ (assume w.l.o.g. that $x_{01}$ is the first non-zero component of $\mb{x_0}$ since $\mb{x_0}\neq 0$). Write $Z_0:=y/x_1-\min\{|z_0|\mb{I}_{|z_0|\neq1/2}-1,1\}$, with $\mb{x}=(x_1,\cdots, x_p)'\in\R^p$.
\vskip 3mm
\noin
\tb{Proposition 3.3} With the same $T$ and $S$ as in Proposition 3.2  under its assumption $1^o)$, further assume that $y$ is symmetrically distributed,
 $\mb{x}$ is spherically distributed, and the distribution of $Z:=(y-\mb{x'}\bs{\beta})/\mb{x'}\mb{v}$ is differentiable near $0$  with density $f_Z$
 at any given $\bs{\beta}\in\R^p$ and $\mb{v}\in\mbs^{p-1}$.
 Then
\begin{itemize}
\item[] (i)\[
\text{IF}((y_0,\mb{x_0}); T^*, F_{(y,\mb{x})})=
\left(\frac{\min\{|z_0|\mb{I}_{|z_0|\neq1/2}-1,1\})}{2f_{Z_0}(0)F_{y}^{-1}(3/4)}, 0,\cdots, 0\right)\in \R^p,
\]
\vskip 3mm
\item[] (ii)\hspace*{20mm} $\displaystyle\gamma^*(T^*, F_{(y,\mb{x})})=\sup_{z_0\in\R}\frac{|\min\{|z_0|\mb{I}_{|z_0|\neq1/2}-1,1\}|}{2f_{Z_0}(0)F_{y}^{-1}(3/4)},$
\end{itemize}
\vskip 3mm
\noindent
\tb{Proof}: (i) Assume, in virtue of equivariance, that $T^*(F)=0$. Then for $\mb{z}=(y_0,\mb{x_0})$ we have 
\be
\text{IF}(\mb{z}; T^*, F)=\lim_{\ep\to 0^+}\frac{T^*(F_{(y,\mb{x})}(\ep,\delta_{\mb{z}}))}{\ep},\label{T*-if-eqn}
\ee
and that
\begin{eqnarray}
T^*(F_{(y,\mb{x})}(\ep,\delta_{\mb{z}}))&=&\arg\!\!\min_{\bs{\beta}\in\R^p} \sup_{\|\mb{v}\|=1}\frac{|T(F_{(\mb{v},\bs{\beta})}(\ep, \delta_{\mb{z}}))|}{S(F_y(\ep,\delta_{y_0}))}\nonumber\\[1ex]
&=&\frac{\arg\!\min_{\bs{\beta}\in\R^p} \sup_{\|\mb{v}\|=1}{|T(F_{(\mb{v},\bs{\beta})}(\ep, \delta_{\mb{z}}))|}}{S(F_y(\ep,\delta_{y_0}))},\label{T^*-expression-1-eqn}
\end{eqnarray}
where $F_{(\mb{v},\bs{\beta})}:=F_{(y-\mb{x}'\bs{\beta}, \mb{x}'\mb{v})}$.
\vskip 3mm
The (L-iv) of Lemma 3.1 can be employed to take care of the denominator of (\ref{T^*-expression-1-eqn}).  In fact,
it tends to $F_y^{-1}(3/4)$
 as $\ep\to0^+$ by the  Lemma 3.1 and the given conditions. We now focus on the numerator of the RHS of (\ref{T^*-expression-1-eqn}).\vskip 3mm

It is readily seen that the distribution of $Z$ is the same for any $\mb{v}\in\mbs^{p-1}$ and a given $\bs{\beta}\in\R^p$  and hence is symmetric about the origin.
By the (L-ii) of Lemma 3.1,
write $T$ in the numerator of the RHS of (\ref{T^*-expression-1-eqn}) for a given $\mb{v}=(v_1,\cdots, v_p)\in\mbs^{p-1}$ ($\mb{x_0'}\mb{v}\neq 0$) and a $\bs{\beta}=(\beta_1,\cdots, \beta_p)\in\R^p$ as 
\[
T(F_{(\mb{v},\bs{\beta})}(\ep, \delta_{\mb{z}}))
=\text{Med}(F_{(\mb{v},\bs{\beta})}(\ep, \delta_{\mb{z}}))=\text{Med}\{A, B, \eta\},
\]
 where $\eta=(y_0-\mb{x_0'}\bs{\beta})/\mb{x_0'}\mb{v}$ 
 and
$A=F_Z^{-1}(1-q(\ep))$ and $B=F_Z^{-1}(q(\ep))$ as defined in Lemma 3.1 and $q(\ep)=1/(2(1-\ep)$. By a direct derivation or standard result on the influence function of the median functional (e.g. Example 3.1 of Huber (1981)), we have
\[
\lim_{\ep \to 0^+}\frac{T(F_{(\mb{v},\bs{\beta})}(\ep, \delta_{\mb{z}}))}{\ep}=
\left\{
\begin{array}{ll}
\frac{-1}{2f_Z(F_Z^{-1}(1/2))},& \text{~if $\eta <F_Z^{-1}(1/2)$}\\[2ex]
0,& \text{~if $\eta =F_Z^{-1}(1/2)$}\\[2ex]
\frac{1}{2f_Z(F_Z^{-1}(1/2))}, & \text{~if $\eta > F_Z^{-1}(1/2)$}\\[1ex]
\end{array}
\right.
\]
Note that by the symmetry of the distribution of $Z$, $F_Z^{-1}(1/2)=0$.
\vskip 3mm

For the consideration of the supremum within the numerator of the RHS of (\ref{T^*-expression-1-eqn}), we should ignore the case $\eta =0$ and just focus on the case $\eta\neq 0$.
Note that the distribution of $Z$ is identical for any $\mb{v}\in\mbs^{p-1}$ and a given $\bs{\beta}\in\R^p$, it is readily seen that if
$\eta\neq 0$, then
\be
\lim_{\ep \to 0^+}\!\!\sup_{\|\mb{v}\|=1}\frac{|T(F_{(\mb{v},\bs{\beta})}(\ep, \delta_{\mb{z}}))|}{\ep}=\sup_{\|\mb{v}\|=1}\lim_{\ep \to 0^+}\frac{|T(F_{(\mb{v},\bs{\beta})}(\ep, \delta_{\mb{z}}))|}{\ep}
=\frac{1}{2f_Z(0)}. \label{prop-2.3.eqn}
\ee
Note that the RHS of (\ref{prop-2.3.eqn}) depends on $\bs{\beta}$ only through the definition of $Z$ and $\eta$.
In order to overall minimize the RHS of (\ref{T^*-expression-1-eqn}),
obviously we have to select $\bs{\beta}$ so that 
$f_Z(0)$ is maximized meanwhile $\eta\neq 0$.
But for any given $\bs{\beta}$ the distribution of $Z$ is symmetric about the origin and its density is maximized at
the origin. Therefore,
 $\bs{\beta}=(\beta_1, 0,\cdots, 0)\in\R^p$ with $\beta_1= \min\{|z_0|\mb{I}_{|z_0|\neq1/2}-1, 1\}$ is obviously one solution.
\vskip 3mm
By the given condition, w.l.o.g., we can select $\mb{v}=(1,0,\cdots, 0)\in\mbs^{p-1}$ in the above discussion and in the definition of $Z$. Then $Z=(y-\mb{x'}\bs{\beta})/(\mb{x'}\mb{v})=Z_0$ and $\eta=z_0-\beta_1\neq 0$.
This, in conjunction with (\ref{T*-if-eqn}) and (\ref{T^*-expression-1-eqn}), 
yields the desired result (i).\vs
(ii) This part is trivial.
\hfill \pend
\vskip 3mm

\noindent
\tb{Remarks 3.4}\vskip 3mm

(I) The influence functions of the P-estimates in MY93
have never been established.
\vskip 3mm
(II) 
Having a bounded influence function or even bounded gross error sensitivity is a very much desirable property for any regression estimating functional. The proposition shows that the deepest projection regression depth  functional $T^*$ possesses this desired property.
\vskip 3mm

(III) The IF of  the deepest regression depth estimating functional in RH99,  has been investigated in VAR00. Where the authors started with elliptical symmetric $(\mb{x}, y)$ but with an appropriate
transformation, the problem is converted to the one with a spherical symmetric $(\mb{x},y)$ for the IF of any regression, scale, affine equivariant functional.  A rather complicated yet bounded IF when $\mb{x}\in\R$ (i.e. $p=1$ here, the simple regression case)  was obtained.
\vskip 3mm

(IV)
The symmetry assumption of the distribution of $y$  could be dropped, then
$F^{-1}_y(3/4)$ in the proposition should be replaced by $F_{|y-c|}^{-1}(1/2)$ with $c=F_y^{-1}(1/2)$.
\hfill\pend
\vskip 3mm
\subsection{Finite sample breakdown point}

Asymptotic breakdown point (ABP) measures the global robustness of a regression estimating functional.
It does not reveal the effect of dimension $p$ on its breakdown point robustness, notwithstanding. In finite sample real practice, there is an alternative to ABP. 
\vskip 3mm

Donoho (1982) and Huber and Donoho (1983) (DH83) introduced the notion of the \emph{finite sample breakdown point} (FSBP) which has become  the most prevailing quantitative measure of global robustness of any location and regression estimators in the finite sample practice.
\medskip

Roughly speaking, the FSBP is the minimum fraction of `bad' (or contaminated) data that the estimator can be affected to an arbitrarily large extent. For example, in the context of estimating the center of a distribution,
the mean has a breakdown point of $1/n$ (or $0\%$), because even one bad observation can change the mean
by an arbitrary amount; in contrast, the median has a breakdown point of $\lfloor(n+1)/2\rfloor/n$ (or $50\%$), where $\lfloor \cdot\rfloor$ is the floor function. For a discussion on general upper and lower bounds of FSBP, see C. M\"{u}ller (2013).
\vskip 3mm

\noindent
\textbf{Definition}
The finite sample \emph{replacement breakdown point} (RBP) of a regression estimator T at the given sample
$Z^{(n)}=\{Z_1,Z_2,\cdots,Z_n\}$, where $Z_i:=(y_i, \mb{x}_i')$, is defined  as
\begin{equation}
\text{RBP}(T,Z^{(n)}) = \min_{1\le m\le n}\bigg\{\frac{m}{n}: \sup_{Z_m^{(n)}}\|T(Z_m^{(n)})- T(Z^{(n)})\| =\infty\bigg\},
\end{equation}
where $Z_m^{(n)}$
denotes an arbitrary contaminated sample by replacing $m$ original sample points in $Z^{(n)}$ with arbitrary points in $\R^{p+1}$.
 Namely, the RBP of an estimator is the minimum replacement fraction which could drive
the estimator beyond any bound. \vskip 3mm

We shall say  $Z^{(n)}$  is\emph{ in general position}
when any $p$ of observations in $Z^{(n)}$ give a unique determination of $\bs{\beta}$.
In other words, any (p-1) dimensional subspace of the space $(y, \mb{x'})$ contains at most p observations of
$Z^{(n)}$.
When the observations come from continuous distributions, the event ($Z^{(n)}$ being in general position) happens with probability one.

\vskip 3mm
\noindent
\textbf{Proposition 3.4} For $T^*$ defined in (\ref{eqn.T*})
with $(T, S)=(\text{Med},~\text{MAD})$  and $Z^{(n)}$ being in general position, 
we have for $ 1\leq p<\lfloor{n/2}\rfloor+2$
\be
\text{RBP}(T^*, Z^{(n)})=\left\{
\begin{array}{ll}
\lfloor (n+1)/2\rfloor\big/n, & \text{if $p=1$,}\\[1ex]
(\lfloor{n}/{2}\rfloor-p+2)\big/n,& \text{if $p>1$,}\\
\end{array}
\right. \label{T*-bp.eqn}
\ee
\vskip 3mm

\noindent
\textbf{Proof:}\vskip 3mm

Note that when $p=1$, the problem becomes an estimation of a location parameter ${\beta_0}$ of $y$ based on minimizing $|\text{Med}_i\{y_i-\beta_0\}|$, and the solution is the median of $\{y_i\}$ which indeed has a RBP given in (\ref{T*-bp.eqn}). In the following, we consider the case $p>1$.
\vskip 3mm
(i) First, \emph{{we show that}} \emph{$m=\lfloor{n}/{2}\rfloor-p+2$ points are enough to breakdown $T^*$}.
 Recall the definition of $T^*(Z^{(n)})$. One has
\begin{eqnarray}
T^*(Z^{(n)})&=&\arg\!\!\min_{\bs{\beta}\in\R^{p+1}}\sup_{\|v\|=1}\frac{\bigg|\text{Med}_{i,\mb{w_i'}\mb{v}\neq 0} \left\{ \frac{y_i-\mb{w_i'}\bs{\beta}}{\mb{w_i'}\mb{v}}\right\}\bigg| }{\text{MAD}_{\substack{1\leq i\leq n}}\{y_i\}}\nonumber  \\[2ex]
&=&\frac{\arg\!\!\min_{\bs{\beta}\in\R^{p+1}}\sup_{\|v\|=1}{\bigg|\text{Med}_{i,\mb{w_i'}\mb{v}\neq 0} \left\{ \frac{y_i-\mb{w_i'}\bs{\beta}}{\mb{w_i'}\mb{v}}\right\}\bigg|}}{\text{MAD}_{\substack{1\leq i\leq n}}\{y_i\} }.
\label{T*-bp-proof.eqn}
\end{eqnarray}

\vskip 3mm
Select $p-1$ points from $Z^{(n)}=\{y_i, \mb{x_i'}\}$. They, together with the origin, form a $(p-1)$-dimensional
subspace (hyperline) $L_h$ in the $(p+1)$-dimensional space of $(y,\mb{x})$.\vskip 3mm
 (Note that since our model contains an intercept term, 
  we assume that the observation ${Z_i}=\mb{0}$ has been deleted from $Z^{(n)}$ 
  for it provides no information on the parameter $\bs{\beta}$).\vskip 3mm
 Construct a non-vertical hyperplane $H$ through $L_h$ (that is, it is not perpendicular to the horizontal hyperplane $y=0$). Let $\bs{\beta}$ be determined by the hyperplane $H$ through $y=\mb{w'}\bs{\beta}$.\vskip 3mm

 We can tilt the hyperplane $H$ so that it approaches its ultimate vertical position. Meanwhile we put all the
 $m$ contaminating points onto this hyperplane $H$ so that it contains no less than $m+(p-1)=\lfloor {n/2}\rfloor+1$
 observations. Call the resulting contaminated sample by $Z^{(n)}_m$.
 Therefore the majority of $(y_i-\mb{w_i'}\bs{\beta})/{\mb{w_i'}\mb{v}}$ now will be zero.
 \vskip 3mm

This implies that $\bs{\beta}$ is the solution for $T^*(Z^{(n)}_m)$ at this contaminated data $Z^{(n)}_m$ since it attains the minimum possible value (zero) on the RHS of (\ref{T*-bp-proof.eqn}). When $H$ approaches  its ultimate vertical position,  $\|\bs{\beta}\|\to \infty$ (for the reasoning, see the proof of Proposition 2.4 of Z18). That is, $m=\lfloor{n}/{2}\rfloor-p+2$ contaminating points are enough to break down $T^*$.  \vskip 3mm

(ii) Second,\emph{ we now show that}\emph{ $m=\lfloor{n}/{2}\rfloor-p+1$ points are not enough to breakdown $T^*$.}
Let $Z^{(n)}_m$ be an arbitrary contaminated sample and $\bs{\beta_c}:=T^*(Z^{(n)}_m)$ and $\bs{\beta_o}=T^*(Z^{(n)})$,  where $Z^{(n)}=\{Z_i\}=\{y_i, \mb{x_i'}\}$ are uncontaminated original points and $\mb{w_i'}=(1, \mb{x_i'})$.
 Assume that $\bs{\beta_c}\neq \bs{\beta_o}$ (Otherwise, we are done). It suffices to show that $\|\bs{\beta_c}-\bs{\beta_o}\|$ is bounded.\vskip 3mm
 Note that since $n-m= \lfloor(n+1)/2\rfloor+p-1$, the denominator of (\ref{T*-bp-proof.eqn}) is the same for contaminated $Z^{(n)}_m$ or original $Z^{(n)}$. We thus focus on its numerator of the RHS of (\ref{T*-bp-proof.eqn}).
~Define
\bee
\delta &=&\frac{1}{2}\inf~\big\{\tau>0; ~\mbox{$\exists$ a $(p-1)$-dimensional subspace  $L$ of ($y=0$) such} \\
& & \mbox{that ${L}^{\tau}$ contains at least $p$ of uncontaminated  $Z_i=(y_i,\mb{x'_i})$ in $Z^{(n)}$}\big\},
\ene
where $L^{\tau}$ is the set of all points $z=(y,\mb{x'})$ that have the distance to $L$ no greater than $\tau$.
Since $Z^{(n)}$ is in general position,  $\delta>0$.\vskip 3mm
Let $H_o$ and $H_c$ be the hyperplanes determined by $y=\mathbf{w}'\boldsymbol{ \beta_o}$ and $y=\mathbf{w}'\boldsymbol{\beta_c}$, respectively, and $M=\max_{i}\{|y_i-\mb{w_i'}\mb{\bs{\beta}}|\}$ for all original $y_i$ and $\mb{x_i}$ in $Z^{(n)}$ with $\mb{w_i'}=(1,\mb{x_i'})$.  Since $\bs{\beta_o}\neq \bs{\beta_c}$, then $H_o\neq H_c$.\vskip 3mm

\textbf{{(A)}} \textbf{Assume that $H_o$ and $H_c$ are not  parallel}.
Denote the vertical projection of the intersection
$H_o\cap H_c$ to the horizontal hyperplane $y=0$ by $L_{vp}(H_o\cap H_c)$, then it is $(p-1)$-dimensional.
By the definition of $\delta$, there are at most $p-1$ of points of ${Z_i}$ within $L_{vp}^{\delta}(H_o\cap H_c)$. Denote the set of all these  possible  $Z_i$ (at most $p-1$) by $S_{cap}$ and $|S_{cap}|=n_{cap}$.
 where ``$|\cdot|$" stands for the counting measure for a set.
 Denote the set of all remaining uncontaminated ${Z_i}$  from the original $\{{Z_i}, i=1,\cdots, n\}$ by $S_r$ and the set of all such $i$ as $I$, then there are at least $n- m-n_{cap}\geq n-\lfloor{n/2}\rfloor=\lfloor{(n+1)/2}\rfloor$ such ${Z}_i$ in $S_r$.
\vskip 3mm

For each $(y_i, \mb{x_i})$ with $i\in I$, construct a two dimensional vertical plane $P_i$ that goes through $(y_i, \mb{x_i})$ and $(y_i+1, \mb{x_i})$ and is perpendicular to  $L_{vp}(H_o\cap H_c)$. Denote the angle formed by $H_o$ and the horizontal line in
$P_i$  by $\alpha_0\in(-\pi/2, \pi/2)$, similarly by $\alpha_c$ for $H_c$ and $P_i$. These are essentially the  angles formed between $H_o$ and $H_c$ with the horizontal hyperplane $y=0$, respectively.  \vskip 3mm
 We see that for $i\in I$ and each  $(y_i, \mb{x_i})$,
$|\mb{w'_i}\bs{\beta_o}|>\delta|\tan(\alpha_o)|$ and $|\mb{w'_i}\bs{\beta_c}|>\delta|\tan(\alpha_c)|$ (see Figure 15 of Rousseeuw and Leroy (1987) (RL87) of a geographical illustration for better understanding, $\mb{x}$ there is $\mb{w}$ here) and $\|\bs{\beta_o}\|=|\tan(\alpha_o)|$ and $\|\bs{\beta_c}\|=|\tan(\alpha_c)|$.
\vskip 3mm

For a given $\mb{v}\in\mbs^{p-1}$ such that $\mb{w'_i}\mb{v}\neq0$ for all $i=1,\cdots, n$.
Write $K_{M}=\min_{i}\{|\mb{w_i'}\mb{v}|\}$ for the given $\mb{v}$ and $K_S=\sup_{i,\mb{v}\in\mbs^{p-1}}\{|\mb{w_i'}\mb{v}|\}$, where $\mb{w_i}=(1,\mb{x_i}')'$ are based on the original uncontaminated $\mb{x_i}$. Then $K_{M}>0$.\vskip 3mm
 Now for each $i\in I$ and the given $\mb{v}$, denote $r^o_i: =(y_i-\mb{w'_i}\bs{\beta_o})/\mb{w'_i}\mb{v}$ and $r^c_i: =(y_i-\mb{w'_i}\bs{\beta_c})/\mb{w'_i}\mb{v}$.
 For the given $\mb{v}$ and any $i\in I$, it follows that (see Figure 15 of RL87)
 \bee
 |r^o_i-r^c_i|&=&\bigg|\frac{\mb{w'_i}\bs{\beta_o}-\mb{w'_i}\bs{\beta_c}}{\mb{w'_i}\mb{v}}\bigg|
 ~~\geq~~ \frac{\delta |\tan(\alpha_o)-\tan(\alpha_c)|}{|\mb{w'_i}\mb{v}|}\\[1ex]
 &\geq& \frac{\delta \big| |\tan(\alpha_o)|-|\tan(\alpha_c)|  \big|}{|\mb{w'_i}\mb{v}|}
 ~~=~~\frac{\delta\big|\|\bs{\beta_o}\|-\|\bs{\beta_c}\|\big|}{|\mb{w'_i}\mb{v}|}\\[1ex]
&\geq& \frac{\delta\big|\|\bs{\beta_o}-\bs{\beta_c}\|-2\|\bs{\beta_o}\|\big|}{|\mb{w'_i}\mb{v}|}
 \ene
 \vskip 3mm

If we assume that $\|\bs{\beta_o}-\bs{\beta_c}\|\geq 2(\|\bs{\beta_o}\|+MK/\delta)$, where $K\geq (K_S+K_M)/2K_{M}$, then by the inequality above we have for $i\in I$ and the given $\mb{v}$
$$
|r^o_i-r^c_i|\geq \frac{\delta\big|\|\bs{\beta_o}-\bs{\beta_c}\|-2\|\bs{\beta_o}\|\big|}{|\mb{w'_i}\mb{v}|}
\geq 2MK/{|\mb{w'_i}\mb{v}|}
$$
which implies that for any $i\in I$ and the given $\mb{v}$,
\[ |r^c_i|\geq |r^o_i-r^c_i|-|r^o_i| \geq  \frac{2MK}{|\mb{w'_i}\mb{v}|}-\frac{ M}{|\mb{w'_i}\mb{v}|}\geq\frac{(2K-1)M}{K_S}\geq \frac{ M}{K_{M}},
\]
which further implies that for the contaminated $(y_i, \mb{x'_i})$ in $Z_m^{(n)}$ and the given $\mb{v}$, we have
\[
\bigg|\text{Med}_{\mb{w'_i}\mb{v}\neq 0 }\left\{ \frac{y_i-\mb{w_i'}\bs{\beta_c}}{\mb{w'_i}\mb{v}} \right\} \bigg|\geq
\frac{M}{K_{M}},
\]
since there are at least $\lfloor(n+1)/2\rfloor$ many $i$ in $I$.\vskip 3mm
 On the other hand, for the given $\mb{v}$, if we
compare all $$\left\{r_i^c\big(\bs{\beta_o}; Z^{(n)}_m\big):=(y_i-\mb{w'_i}\bs{\beta_o})/(\mb{w'_i}\mb{v})\right\},~\text{ where $(y_i, \mb{x'}_i)$ is from $Z_m^{(n)}$},$$
with all $$\left\{r^o_i\big(\bs{\beta_o}; Z^{(n)}\big) :=(y_i-\mb{w'_i}\bs{\beta_o})/(\mb{w'_i}\mb{v})\right\},~\text{ where $(y_i, \mb{x'}_i)$ is from $Z^{(n)}$},$$
it is readily seen that there are at least $N$ terms are the same, where $N= n_{cap}+|S_r|=n-m$ ($n_{cap}$ original points in $S_{cap}$ plus $|S_r|$ original points in $S_r$).
Therefore,
among all $\{\big|r^c_i\big(\bs{\beta_o}; Z^{(n)}_m\big)\big|\}$, there are at least $n-m\geq (p-1)+\lfloor(n+1)/2\rfloor$ terms each of which is no greater than $M/K_{M}$ since for all $i$, $\big|r^o_i\big(\bs{\beta_o}; Z^{(n)}\big)\big|\leq M/K_{M}$. That is, for $(y_i, \mb{x'}_i)$ from $Z_m^{(n)}$ and the given $\mb{v}$
\be
\bigg|\text{Med}_{i,\mb{w'_i}\mb{v}\neq 0 }\left\{ \frac{y_i-\mb{w_i'}\bs{\beta_o}}{\mb{w'_i}\mb{v}} \right\} \bigg|
\leq M/K_{M}. \label{bp-upper-eqn}
\ee

\vskip 3mm
Assume that $\mb{v}$ is the direction at which $\bs{\beta_c}$ attains the minimum of  the numerator of the RHS of (\ref{T*-bp-proof.eqn}). That is, for $(y_i, \mb{x'}_i)$ from $Z_m^{(n)}$
\[
\inf_{\bs{\beta}\in\R^{p+1}}\sup_{\|\mb{v}\|=1}\bigg|\text{Med}_{i,\mb{w'_i}\mb{v}\neq 0}\left\{\frac{y_i-\mb{w'_i}\bs{\beta}}{\mb{w'_i}\mb{v}}\right\}\bigg|
=\bigg|\text{Med}_{i, \mb{w'_i}\mb{v}\neq 0}\left\{\frac{y_i-\mb{w'_i}\bs{\beta_c}}{\mb{w'_i}\mb{v}}\right\}\bigg|,
\]
 Hence it follows that for $(y_i, \mb{x'}_i)$ from $Z_m^{(n)}$ and the $\mb{v}$
\[
\bigg|\text{Med}_{\mb{w'_i}\mb{v}\neq 0}\left\{ \frac{y_i-\mb{w_i'}\bs{\beta_c}}{\mb{w'_i}\mb{v}} \right\} \bigg|\leq \bigg|\text{Med}_{\mb{w'_i}\mb{v}\neq 0}\left\{ \frac{y_i-\mb{w_i'}\bs{\beta_o}}{\mb{w'_i}\mb{v}} \right\} \bigg|\leq
\frac{M}{K_{M}},
\]
The first inequality follows from the definition of $\bs{\beta_c}$ and $v$, the second one follows from
the inequality  (\ref{bp-upper-eqn}) established above. Now we reach a contradiction. \vskip 3mm
Therefore, $\|\bs{\beta_o}-\bs{\beta_c}\|< 2(\|\bs{\beta_o}\|+MK/\delta)$
and thus $\|\bs{\beta_o}-\bs{\beta_c}\|$ is bounded. That is, $m$ contaminating points are not enough to breakdown $T^*$. \vskip 3mm

\textbf{{(B)}} \textbf{Assume that $H_o$ and $H_c$ are parallel}. That is,
$\bs{\beta_c}=\rho\bs{\beta_o}$. If $\rho$ is finite, then $\|\bs{\beta_c}-\bs{\beta_o}\|$ is automatically bounded. We are done. Now consider the case that $|\rho| \to\infty$, that is, $|\rho|$ can be arbitrarily large.\vskip 3mm

\tb{(B1)} \textbf{Assume that $H_o$ is not parallel to $y=0$}.\vskip 3mm
 The proof is very similar to part \tb{(A)}.
Denote the intersection of $H_c$ and the horizontal hyperplane $y=0$: $H_c\cap\{y=0\}$ by $L_c$. Then
$L^\delta_c$ contains at most $p-1$ uncontaminated points from $\{Z^{(n)}\}$. Denote the set of all the remaining uncontaminated points in $\{Z^{(n)}\}$ as $S_r$. Hence $|S_r|\geq n-m-(p-1)\geq \lfloor(n+1/2\rfloor$.
Denote again by $I$ the set of all $i$ such that $Z_i\in S_r$. Again let the angle between $H_c$ and $y=0$ be $\alpha_c$, then it is seen that $\|\bs{\beta_c}\|=|\tan(\alpha_c)|$ and $|\mb{w'_i}\mb{\bs{\beta_c}}|> \delta
|\tan(\alpha_c)|$ for any $i\in I$.
\vskip 3mm

Assume that $\mb{v}_c$ is one unit vector at which $\bs{\beta_c}$ attains the $\inf$ of the numerator of the HRS of (\ref{T*-bp-proof.eqn}). Define $K_{M}=\min_{i}\{\mb{w_i'}\mb{v_c}\}$, then $K_M>0$. Write 
$$r^c_i=(y_i-\mb{w_i'}\mb{\bs{\beta_c}})/(\mb{w_i'}\mb{v_c}),$$ for all $\mb{x}_i$ (and hence $\mb{w}_i$) from $Z^{(n)}_m=(y_i,\mb{x_i})$.  Write $M_y=\max_{i}|y_i|$. It follows that for $i\in I$
\[
\big|r^c_i\big|\geq \big||\mb{w'_i}\bs{\beta_c}|-|y_i|\big|/K_S\geq |~\delta|\tan(\alpha_c)|-M_y|/K_S.
\]
Since $|S_r|\geq \lfloor(n+1/2\rfloor$, then for all $(y_i, \mb{x'}_i)$ (and hence $\mb{w}_i$) from $Z^{(n)}_m=(y_i,\mb{x'_i})$
$$\bigg|\text{Med}_{i }\left\{ \frac{y_i-\mb{w_i'}\bs{\beta_c}}{\mb{w'_i}\mb{v_c}} \right\} \bigg|\geq
|~\delta|\tan(\alpha_c)|-M_y|/K_S.
$$
Now introduce $r_i^c\big(\bs{\beta_o}; Z^{(n)}_m\big)$ and $r_i^o\big(\bs{\beta_o}; Z^{(n)}\big)$ as in the proof of part \textbf{(A)}.  Therefore,
among all $\{\big|r^c_i\big(\bs{\beta_o}; Z^{(n)}_m\big)\big|\}$, there are at least $n-m\geq (p-1)+\lfloor(n+1)/2\rfloor$ terms  each of which  is no greater than $M/K_{M}$ since for all $i$, $\big|r^o_i\big(\bs{\beta_o}; Z^{(n)}\big)\big|\leq M/K_{M}$. That is, for $(y_i, \mb{x'}_i)$ from $Z_m^{(n)}$ and the given $\mb{v_c}$
\be
\bigg|\text{Med}_{i}\left\{ \frac{y_i-\mb{w_i'}\bs{\beta_o}}{\mb{w'_i}\mb{v_c}} \right\} \bigg|
\leq M/K_{M}.\label{bp-upper-2-eqn}
\ee
\vskip 3mm

On the other hand,
it is not difficult to see that for $(y_i, \mb{x'}_i)$ from $Z_m^{(n)}$
\bee
\inf_{\bs{\beta}\in\R^{p+1}}\sup_{\|\mb{v}\|=1}\bigg|\text{Med}_{i,\mb{w'_i}\mb{v}\neq 0}\left\{\frac{y_i-\mb{w'_i}\bs{\beta}}{\mb{w'_i}\mb{v}}\right\}\bigg|
&=&\bigg|\text{Med}_{i}\left\{\frac{y_i-\mb{w'_i}\bs{\beta_c}}{\mb{w'_i}\mb{v_c}}\right\}\bigg|\\[2ex]
&\leq& \bigg|\text{Med}_{i}\left\{ \frac{y_i-\mb{w_i'}\bs{\beta_o}}{\mb{w'_i}\mb{v_c}} \right\} \bigg|~~\leq~~
\frac{M}{K_{M}},\\[0ex]
\ene
where the first inequality follows directly from the definitions of $\bs{\beta_c}$ and $\mb{v_c}$ and the second
one directly from  (\ref{bp-upper-2-eqn}).\vskip 3mm
If $|\rho|$ could be arbitrarily large, then since $\delta|\tan(\alpha_c)|-M_y=\delta|\rho|\|\bs{\beta_o}\|-M_y$ could be arbitrarily large, so that  $|~\delta|\tan(\alpha_c)|-M_y|/K_S> M/K_M$,
which leads to a contradiction.  Hence $\|\bs{\beta_o}-\bs{\beta_c}\|$ is bounded. It means that $m$ contaminating points are not enough to breakdown $T^*$.

 \vskip 3mm

\tb{(B2)} \textbf{Assume that $H_o$ is parallel to $y=0$}.
Then, it means that $\bs{\beta_c}=\rho\bs{\beta_o}=(\rho\beta_{o1},0,\cdots, 0)$. Assume that $\beta_{o1}\neq 0$. Otherwise, we are done.  Now we can repeat the argument above since $n-m\geq (p-1)+\lfloor(n+1)/2\rfloor$.  On the one hand  we can show that for all $(y_i, \mb{x'}_i)$ from $Z_m^{(n)}$
\bee
\bigg|\text{Med}_{i}\left\{\frac{y_i-\mb{w'_i}\bs{\beta_c}}{\mb{w'_i}\mb{v_c}}\right\}\bigg| 
~~\leq~~ \bigg|\text{Med}_{i}\left\{ \frac{y_i-\mb{w_i'}\bs{\beta_o}}{\mb{w'_i}\mb{v_c}} \right\} \bigg|~~\leq~~
\frac{M}{K_{M}}, &&\\[0ex]
\ene
where, $\bs{\beta_o}$ and $\bs{\beta_c}$, $M$ and $K_M$ and $\mb{v_c}$ are defined as before.
\vskip 3mm
On the other hand, we have for all $(y_i, \mb{x'}_i)$  from $Z^{(n)}_m=(y_i,\mb{x'_i})$
$$\bigg|\text{Med}_{i }\left\{ \frac{y_i-\mb{w_i'}\bs{\beta_c}}{\mb{w'_i}\mb{v_c}} \right\} \bigg|\geq
\big|~|\rho\beta_{o1}|-M_y\big|\big/K_S, $$
where $K_S$ and$M_y$ is defined as before.
\vskip 3mm

Again if $|\rho|$ could be arbitrarily large, then since $|\rho\beta_{o1}|-M_y$ could be arbitrarily large so that  $\big||\rho\beta_{o1}|-M_y\big|/K_S> M/K_M$
yields a contradiction. Hence $\|\bs{\beta_o}-\bs{\beta_c}\|$ is bounded.
That is, $m$ contaminating points are not enough to breakdown $T^*$.
 \hfill \pend

\vskip 3mm

\noindent
\textbf{Remarks 3.5}\vskip 3mm
(I) MY93 also discussed the FSBP of their P-estimates,  the RBP of the P-estimates has never been established, nevertheless. MY93 established an upper bound for the norm of their P-estimates which holds true with some probability that could be very close to one by taking sufficiently large number of subsamples in the computation of their P-estimates.  Although P-estimates are defined differently from $T^*$ here, the idea of the proof above, however, seems applicable to the P-estimates to obtain a concrete (and with probability one) RBP.
\vs
(II) The main idea of the proof above was adapted from the proof of the RBP of the LMS in Rousseeuw (1984). 
The latter, however,  only addressed part \tb{(A)}, and part \tb{(B)} was overlooked, where it was assumed implicitly that $H_c\cap H_o\neq \emptyset$. The same assumption was made in the proof of the RBP of the LTS
(page 132 of RL87). One may ask how often in practice
 $H_c\cap H_o=\emptyset$? The argument seems reasonable at first.  However, one cannot afford to miss any conceivable contamination case when establishing RBP.
  \vskip 3mm
(III) Although $T^*_n$ possess a very high RBP (the same as that of LMS), it is still not the best possible RBP for any regression equivariant estimator. For the latter, it is $(\lfloor\frac{n-p}{2}\rfloor+1)/n$ (see page 125 of RL87). To attain the upper bound of RBP, one can modify the $T^*_n$ so that its RBP attains the upper bound.
Indeed, there are several variants of the $T^*_n$ below. \vskip 3mm
First, in the definition of $T^*_n$, consider the
median of all $\big\{|\frac{y_i-\mb{w'_i}\bs{\beta}}{\mb{w'_i}\mb{v}}|\big\}$. That is, consider the median of the absolute values instead of the absolute value of the median. Call the resulting estimator $\text{T1}^*_n$. Second, replace the median of the absolute values by the $h$th ordered absolute values. If $h=\lfloor n/2\rfloor+1$,
call the resulting estimator $\text{T2}^*_n$. If  $h=\lfloor n/2\rfloor+\lfloor(p+1)/2\rfloor$, call the resulting estimator $\text{T3}^*_n$. 
One can show that the RBP of $\text{T1}^*_n$ or $\text{T2}^*_n$ is the same as $T^*_n$ but that of
$\text{T3}^*_n$ attains the upper bound. Thirdly, other variants include replacing the $h$th ordered absolute values with the sum of first $h$th ordered absolute values, then the resulting estimators have the same RBP
of $\text{T2}^*_n$ and $\text{T3}^*_n$, respectively, corresponding to the two choices of $h$:  $h=\lfloor n/2\rfloor+1$ or $h=\lfloor n/2\rfloor+\lfloor(p+1)/2\rfloor$.

\vskip 3mm
(IV) To the best of knowledge of this author,
 the RBP of $T^*_{RD}$ (the deepest regression estimator defined in RH99),  has not yet been established \emph{explicitly}.
\vskip 3mm
(V) The RBP result is established under the assumption that $Z^{(n)}$ is in general position. In more general cases, one can use a number $c(Z^{(n)})$ (which is the maximum number of observations from $Z^{(n)}$ contained in any $(p-1)$ dimensional subspace) to replace $p$ in the derivation of the final RBP result.
\hfill \pend
\vskip 3mm

\section{Computation, robustness illustration, and simulation}
\vs
\subsection{Computation}
The deepest projection regression depth estimator $T_n^*$ faces a common problem for any estimators with high breakdown point robustness. That is, it is very challenging to compute them
 in practice while enjoying the best possible ABP.\vskip 3mm
 \noin
Exact computation of $T^*_n$ is certainly difficult (it involves two layers of optimizations (minimization of the maximized unfitness), if not impossible. But one can at least compute $T^*_n$ approximately. Here sub-sampling schemes and the MCMC technique could be employed in the optimization process, as done in Shao and Zuo (2017) for halfspace depth in $\R^d$.
\vs
\noin
The rough idea is as follows.
Randomly select $N_{\bs{\beta}}$ of  $\bs{\beta}$'s over a very wide range in parameter space $\R^p$, calculate all UF$(\bs{\beta},F^n_Z)$. 
 Sort the latter and select $p+1$ $\bs{\beta}$'s with smallest unfitness. Over the simplex formed by these $p+1$ $\bs{\beta}$ points (in parameter space), search for the point $\bs{\beta}$ with the smallest unfitness (equivalent the deepest regression line or hyperplane).
\vs
\noin
In the above process, we have implicitly take the advantage of the property of PRD$(\bs{\beta};F_Z)$ or UF$(\bs{\beta};F_Z)$.
That is, PRD$(\bs{\beta};F_Z)$ satisfies the property (P3) of Z18 (monotonicity relative to the deepest point). Therefore the depth region
of $\bs{\beta}$ (the set of all $\bs{\beta}$'s with its depth no less than a fixed value)  is convex and nested. Hence, the deepest point(s) must lie in the convex simplex formed by the $p+1$ $\bs{\beta}$ points. When there is more than one deepest point, we can take the average of them, the resulting point will possess the maximum depth.
\vs
The following is an {\textbf{approximate algorithm} (AA)} for the computation of $T^*_n$.
\begin{itemize}
\item[] (A) Randomly select a set of points $\bs{\beta_j}\in\R^p$ over a very wide range of region, $j=1,\cdots, N_{\bs{\beta}}$, where $N_{\bs{\beta}}$ is a tuning parameter of the  total number of the random points.
\item[] (B) For each $\bs{\beta_j}$, randomly select a set of unit directions $\mb{v_k}\in\mbs^{p-1}, k=1,\cdots, N_{\mb{v}}$. $N_{\mb{v}}$ is another tuning parameter. Compute the approximate unfitness of $\bs{\beta_j}$ w.r.t.  $\{ Z^j_{ik} =(y_i-\mb{w_i'}\bs{\beta_j})\big/ (\mb{w_i'}\mb{v_k})\}$ for a fixed $j$, and all $i$ and $k$, where, $i=1,\cdots, n$, $k=1,\cdots, N_{\mb{v}}$\vskip 3mm

\item[] (C) Order the $\bs{\beta}_j$'s according to their depth (or equivalently unfitness) and select the deepest $p+1$ $\bs{\beta}_j$'s.
Search over the closed convex hull formed by these $p+1$ points via common optimization algorithms (e.g. the downhill simplex method, or the MCMC technique) to get the final deepest $\bs{\beta}$ or our approximate $T^*_n$.
\item[] (D) To mitigate the effect of randomness, repeat the steps above (many times) so that the one of $T^*_n$ with the maximum updated projection regression depth is adopted.
\end{itemize}
\vs
\tb{Remarks 4.1}:
 \begin{itemize}
 \item[(I)]  The candidate  (random point) $\bs{\beta}$ can be produced by randomly selecting $p$ points from $Z^{(n)}=\{ (\mb{x_i}, y_i), i=1,\cdots, n\}$
 which (by the general position) determine a unique hyperplane $y=\mb{w'}\bs{\beta}$ containing all  $p$ points.\vs

\item[(II)] If Med and MAD are used for the $(T, S)$, then, the random directions could be selected among those which are perpendicular to the hyperplanes formed by $p$ points  from $Z^{(n)}$. 
\item[(III)] For a better approximation of depth (unfitness) of $\bs{\beta}_j$, tuning (increasing) $N_{\mb{v}}$. For a better approximation of $T^*_n$, tuning
$N_{\bs{\beta}}$.
Continue iterations until it satisfies a stopping rule (e.g. the difference between consecutive depths is less than a cutoff value).
\item[(IV)] The overall worst case time complexity of the algorithm is: step (A)+(B): $O(N_{\mb{v}}N_{\bs{\beta}}n)$, where the linear method is employed to compute the univariate median; step (C):  $O(N_{\bs{\beta}}\log(N_{\bs{\beta}})+N_{\mb{v}}N_{\bs{\beta}}n)$, where over the closed convex hull, step (A) and (B) are assumed to be repeated; step(D) $O(R(N_{\mb{v}}N_{\bs{\beta}}n+N_{\bs{\beta}}\log(N_{\bs{\beta}})))$, where $R$ is the number of replications. The overall cost of the
    algorithm is $O(RN_{\bs{\beta}}(N_{\mb{v}}n+\log(N_{\bs{\beta}})))$.
\item[(V)] Theoetically speaking, the AA is suitable for any $p$. But in high dimensions, the $N_{\bs{\beta}}$ and $N_{\mb{v}}$ should be dependent on $p$ to get better approximation. A larger $N_{\bs{\beta}}$ is more important than a large $N_{\mb{v}}$ since a rough approximation of UF$(\bs{\beta}; F^n_Z)$ is allowed as long as the first (p+1) deepest $\bs{\beta}$'s are correctly identified or the most importantly the convex hull formed contains the 
    deepest point.
      To guarantee the latter, in practice $p$ is limited, say $2\leq p <5$ for the AA unless tuning parameters are chosen dependent on $p$.
\end{itemize}

\subsection{Robustness illustration}
With the approximate algorithm above, we are now in a position
 to better appreciate the outstanding breakdown robustness of the deepest projection depth estimator $T^*_{PRD}$. We illustrate below the performance of the regression lines of the classical least squares, the $T^*_{RD}$ of RH99, and the $T^*_{PRD}$ w.r.t. contamination in a data set.\vs
 \noin
  \tb{Example 4.2.1.}: A small data set (given in Table 9 of RL87) (only for illustration purpose). The original data set contains nine bivariate points, but one point (0,0) provides no information for the regression and therefore is deleted, yielding an eight-point data set.
\vs
\bec
\begin{figure}[h]
    \centering
    \vspace*{-5mm}
\includegraphics[width=\textwidth]{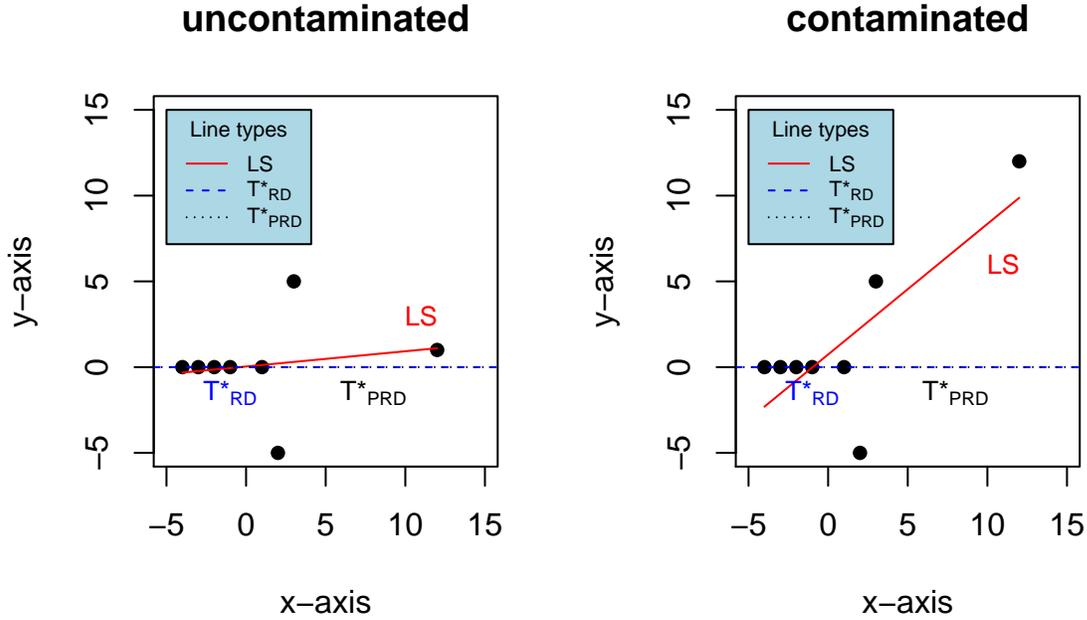}
    \caption{ Three regression lines  for data without or with contamination (red solid line for LS, blue dashed line for $T^*_{RD}$ and black dotted line for $T^*_{PRD}$). Left: Original eight-point data set. $T^*_{RD}$ and $T^*_{PRD}$ are identical.  Right: Contaminated data set with one original point moved form $(12,1 )$ to $(12,12)$, leading to a drastically change in
    the LS line while both $T^*_{RD}$ and $T^*_{PRD}$ are unchanged and resist the contamination.}\vspace*{-0mm}
\end{figure}
\vspace*{-10mm}
\enc
Regression lines given by the three approaches are plotted w.r.t. the original data versus (i) $12.5\%$  contaminated data set (one data point is contaminated) in Figure 1 left and right  and versus (ii)  $37.5\% $ contaminated data set (three points out of eight are contaminated) in Figure 2 left and right, respectively.\vs

Inspecting  Figure 1, reveals that (i) for the original data, the least squares line is affected by the point with large $x$-coordinate (an outlier in the x-direction, or a leverage point). It is drawn by this leverage point, whereas  both deepest regression depth lines resist against the leverage point and capture the horizontal line $y=0$,  (ii) When the leverage point is moved upward to $(12,12)$, then
 the entire least squares line is attracted by this movement and moved upward (which means that a single point can ruin the LS line), whereas  both deepest regression depth lines are resistant to this single point contamination.
\vs
  Figure 2, on the other hand, reveals that (i) for the uncontaminated data, the situation is the same as in Figure 1 left, and (ii) for the contaminated data (three points are contaminated),
 the least squares line again is affected by the leverage point as well as the contaminated points, but not too much from the latter (since the x and y coordinates of the contaminated points are moderate), the deepest line of $T^*_{PRD}$ is affected by the contamination but still informative and useful, whereas the one from $T^*_{RD}$ is useless (breaks down as expected due to more than $1/3$ of contamination). Note that the RD of this vertical line is $4/8$ while there are other lines that have this depth. To deal with the non-uniqueness problem while in order to have the affine equivariance of the final deepest regression line, one can take an average of lines with the maximum depth. But the resulting line will still have a unbounded slope, hence is useless.

\bec
\begin{figure}[h]
    \centering
    \vspace*{-5mm}
\includegraphics[width=\textwidth]{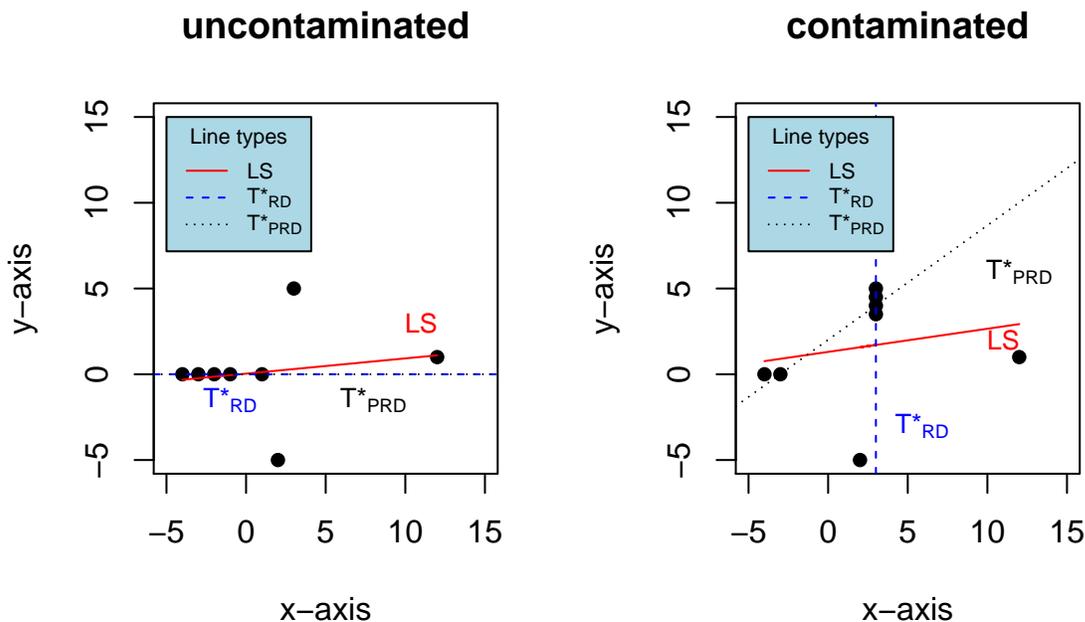}
\vspace*{-10mm}
    \caption{{ Three regression lines  for data without or with contamination (red solid line for LS, blue dashed line for $T^*_{RD}$ and black dotted line for $T^*_{PRD}$). Left: Original eight-point data set, $T^*_{RD}$ and $T^*_{PRD}$ are identical. Right: Contaminated data set with three original points moved to the points with 3 as their x-coordinates, $T^*_{RD}$ breaks down.}}
\end{figure}
\vspace*{-0mm}
\enc
\vs
\vs
\bec
\begin{figure}[h!]
\vspace*{-15mm}
\includegraphics[width=\textwidth]{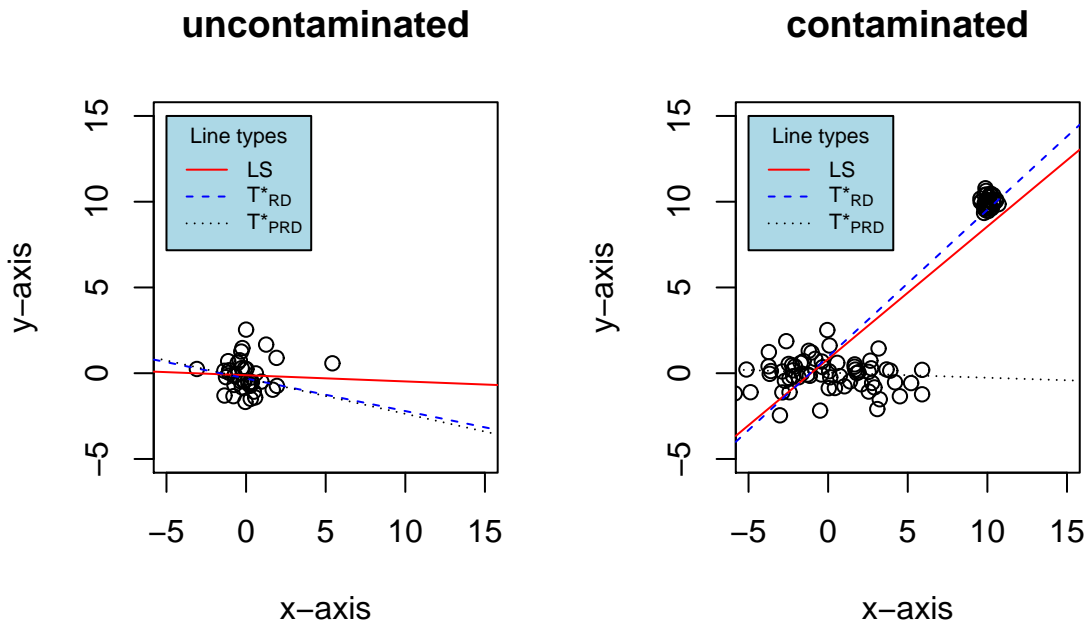}
\vspace*{-10mm}
 \caption{{ Three regression lines for data  without or with contamination (red solid line for LS, blue dashed line for $T^*_{RD}$ and black dotted line for $T^*_{PRD}$). Left: Original 100 normal points, lines from LS, $T^*_{RD}$ and $T^*_{PRD}$ are very similar. Right: Contaminated data set with $34\%$ points contaminated, both LS and $T^*_{RD}$ ``break down" while $T_{PRD}$ resists the contamination and is still useful.}}
 \label{fig:example-4-2}
\end{figure}
\enc
\vspace*{-20mm}
\noin
  \tb{Example 4.2.2.}: We generate a bivariate normal data set with size $100$ and $\mu$ and $\Sigma$  are $$ \mu=\begin{pmatrix}0\\0\end{pmatrix},~~~~\Sigma=\begin{pmatrix}1&-0.8\\-0.8 &1\end{pmatrix};~~~~~~\mu_1=\begin{pmatrix}10\\10\end{pmatrix},~~~~ \Sigma_1=\begin{pmatrix}0.1&0\\0&0.1\end{pmatrix}. $$ Then we consider a $34\%$ replacement normal points contamination with $\mu_1$ and  $\Sigma_1$.   The performance of the three lines is displayed in Figure (\ref{fig:example-4-2}).
\vs

We compute the three lines w.r.t. un-contaminated data. The three (slope, intercept) lines are (-0.11788718, -0.03614133), (-0.3066041, -0.1899350), and (-0.2943452 -0.2073059) for LS, $T_{RD}$, and $T_{PRD}$, respectively. They do not differ very much as shown in the left side of Figure \ref{fig:example-4-2}, or all three seem to be useful.\vs
 On the other hand, we also compute the three lines w.r.t a $34\%$ replacement contamination. The
three lines are (0.8283450, 0.7727246), (0.9657038, 0.8559868 ), and (0.03350186, -0.02969263) for LS, $T_{RD}$, and $T_{PRD}$, respectively.
They differ very much as shown in the right side of Figure \ref{fig:example-4-2}. Both LS and $T_{RD}$ lines break down (attached to the cloud of contamination) whereas $T_{PRD}$ can resist the $34\%$ contamination (in fact up to $50\%$) and continue to provide a useful regression line.
\subsection{Finite-sample relative efficiency}
Robustness does not  work in tandem with efficiency. $T^*_{PRD}$ (or $T^*_n$ in the empirical case) has the best possible ABP while it has to pay a price of a relatively low efficiency. Its efficiency, however, could be improved (as shown below) by replacing, the univariate median,  the chief source of low efficiency, with a much more efficient depth trimmed or weighted mean (Zuo (2006), Zuo, Cui and He (2004) (ZCH04)) meanwhile keeping it as robust as before, just as its location counterpart, the projection median, does (Zuo (2003)).\vs
On the other hand, the deepest regression line in RH99 ($T^*_{RD}$) has no such freedom to improve its low efficiency since it is fixed and unlike $T^*_{PRD}$, which represents a class of
functionals (estimators) with the different choices of univariate functionals $T$ (used in $T^*_{PRD}$) that can be highly efficient yet as robust as the univariate median.\vs

In the following we investigate via simulation the finite-sample relative efficiency of the deepest lines $T^*_{RD}$ and $T^*_{PRD}$ w.r.t. the classical least squares line. We generate $N_R$ samples from the simple linear regression model: $y_i=\beta_0+\beta_1 x_i+e_i, i=1,2,\cdots, n,$
with different sizes $n$ (see Table 1), where $e_i \sim N(0,\sigma^2)$. In light of the regression equivariance, we can assume w.l.o.g. that $\bs{\beta}=(\beta_0, \beta_1)'=(0,0)'$. We generate $x_i$ from standard normal and $t(2)$ independently with $y_i$, which are $ N(0, 1)$ points.
The relative efficiency of the slope and intercept of the lines $T^*_{RD}$ and $T^*_{PRD}$ w.r.t. those of the least squares line are listed in Table 1
with various $n$, where $T$ in the definition of $T^*_{PRD}$ is the sample median.\vs

Inspecting the Table 1 reveals that
(i)
 for Gaussian $x_i$'s the intercept of $T^*_{RD}$ is slightly more efficient than that of  $T^*_{PRD}$ when $n\geq 20$, while the slope of $T^*_{PRD}$ is more efficient than that of $T^*_{RD}$ uniformly for all $n$ , whereas for $t(2)$ $x_i$'s, $T^*_{PRD}$ is more efficient than $T^*_{RD}$ both in slope and intercept uniformly for all $n$;
(ii) the efficiency of the deepest regression lines differs when the $x_i$ are generated from different distributions;
(iii) slopes have higher efficiency for Gaussian $x_i's$ than for $t(2)$ $x_i's$; and (iv) for Gaussian $x_i's$ slopes have higher efficiency than intercept for $n>10$, this relationship is reversed for $t(2)$ $x_i's$ for all $n$.
\begin{table}[h!]
\centering
 Relative efficiency (based on $N_R=5,000$ replications) of the deepest lines $T^*_{RD}$  and $T^*_{PRD}$ compared to the least squares line when the $x_i$ are from Gaussian or $t$ distributions
\begin{tabular}{ccccc}
\hline\\[1ex]
&~~~~~~~~\underline{~~Gaussian $x_i$~~}& & ~~~~\underline{~~t(2) $x_i$~~}&\\[2ex]
n &~~~~~~~~ slope~~~&  intercept &~~ slope & intercept\\[1ex]
 & ~~$(T^*_{PRD}$;~$T^*_{RD})$ &$(T^*_{PRD}$;~$T^*_{RD})$&$(T^*_{PRD}$;~$T^*_{RD})$ &$(T^*_{PRD}$;~$T^*_{RD})$\\[1ex]
10&(0.6990;~0.6341)&(0.7044;~0.6696)&(0.6060;~0.5404)&(0.6950;~0.6789)\\
20&(0.7267;~0.7156)&(0.6884;~0.6941)&(0.6217;~0.6048)&(0.7111;~0.7020)\\
40& (0.7513;~0.7321)&(0.7057;~0.7124)&(0.6154;~0.5967)&(0.7354;~0.7142)\\
80&(0.7606;~0.7514)&(0.7042;~0.7107)&(0.5876;~0.5759)&(0.7285;~0.7063)\\
100 &(0.7471;~0.7385)&(0.7024;~0.7114)&(0.5837;~0.5685)&(0.7286;~0.7155)\\
\hline
\end{tabular}
\caption{
Median used for $T$ in $T^*_{PRD}$}
\label{table-1}
\end{table}

\begin{table}[h!]
\centering
 Relative efficiency (based on $N_R=1,000$ replications) of the deepest line $T^*_{RD}$  and $T^*_{PRD}$ compared to the least squares line when the $x_i$ are from Gaussian or $t$ distributions   \\[2ex]
\begin{tabular}{ccccc}
\hline\\[1ex]
&~~~~~~~~\underline{~~Gaussian $x_i$~~}& & ~~~~\underline{~~t(2) $x_i$~~}&\\[2ex]
n &~~~~~~~~ slope~~~&  intercept &~~ slope & intercept\\[1ex]
 & ~~$(T^*_{PRD}$;~$T^*_{RD})$ &$(T^*_{PRD}$;~$T^*_{RD})$&$(T^*_{PRD}$; ~$T^*_{RD})$ &$(T^*_{PRD}$;~$T^*_{RD})$\\[1ex]
10&(0.6496;~0.6065)&(0.6780;~0.6676)&(0.6561;~0.5874)&(0.7070;~0.6930)\\
20& (0.7196;~0.6982) &(0.7346;~0.7316)&(0.5813;~0.5582)&(0.7311;~0.7010)\\
40&(0.7581;~0.7289)&(0.7784;~0.7655)&(0.5972;~0.5722)&(0.7207;~0.6869)\\
80&(0.7816;~0.7482)&(0.7177;~0.7054)&(0.5992;~0.5686)&(0.7323;~0.7098)\\
100 &(0.7505;~0.7462)&(0.7166;~0.7138)&(0.6201;~0.5978)&(0.7065;~0.6931)\\
\hline
\end{tabular}
\caption{
PWM (see ZCH04) used for $T$ in $T^*_{PRD}$}
\label{table-2}
\end{table}
\vs
The efficiency of the slope and intercept of the line $T^*_{PRD}$ could be improved by replacing median employed in the definition of  $T^*_{PRD}$
with a more efficient projection depth weighted mean (PWM) yet have the same level of robustness as the median, see Zuo (2003) and Zuo, Cui and He (2004) (ZCH04), and Wu and Zuo (2009):
$$
\mbox{PWM}(x^n)=\frac{\sum_{i=1}^n w(PD_n(x_i))x_i}{\sum_{i=1}^n w(PD_n(x_i))},
$$
where $w(r) = I (r < c)\big(\exp(-k(1-r/c)^2)-\exp(-k)\big)/(1-\exp(-k))+I (r \geq c)$,
$PD_n(x_i) = 1/(1+|x_i -\mbox{Med}(x^n)|/\mbox{MAD}(x^n))$ and  $x^n=\{x_1,\cdots, x_n\}$ with $x_i \in \R$.
For discussions of weight function $w$ and parameters $k$ and $c$, see ZCH04.
\vs
 Generally speaking, tuning $c$ to render it smaller to get higher efficiency from
PWM. The same is true  for parameter $k$. Namely, keeping the number of inner points as large as possible to gain higher efficiency and down-weighting outliers slower to gain higher efficiency.  In our simulation, we set $k=3$ and $c=3.5$. Other parameters that could be tuned include $N_v$ and $N_{\beta}$ (see Section 4.1). In our simulation, we set $N_v=100+2*n$, where $100$ are random directions and $2n$ directions (they are $\mb{v_i}\pm (10^{-10}, 0)'$, where $\mb{w}_i'\mb{v}_i=0$, see the RHS of (\ref{T*-bp-proof.eqn}), $i=1,\cdots,n$) are strategically chosen. $N_{\beta}$ is increasing with $n$ but no greater than $n(n-1)/2$. With these parameters, the results from $T^*_{RD}$ and $T^*_{PRD}$ are listed in Table 2.
\vs
Inspecting Table 2 reveals that
(i) with the PWM employed in the definition of $T^*_{PRD}$, $T^*_{PRD}$ becomes more efficient than $T^*_{RD}$ both in slope and intercept uniformly for all $n$ both for Gaussian $x_i$ and $t(2)$ $x_i$
(note that by tuning the parameters, one can even gets higher efficiency for $T^*_{PRD}$).
(ii) The efficiency of the deepest lines depends on the distribution of $x_i$.
(iii) The efficient of the intercept is higher than that of slope for $t(2)$ $x_i$'s. This is no longer true for Gaussian $x_i$'s and when $n>40$.
\section{Discussions and concluding remarks}

This article investigates the  robustness property of the
deepest projection regression depth functional $T^*_{PRD}$. $T^*_{PRD}$ is closely related to (but different from) the P-estimates in MY93. In fact, it is the modification of the latter, to achieve the scale invariance of the induced depth function and scale equivariance of $T^*_{PRD}$.\vs
 Like MY93 for the P-estimates, an upper bound for the maximum bias of $T^*_{PRD}$ is established, which covers Theorems 3.4, 3.5, and  4.1 of MY93.  In contrast to MY93 for their P-estimates, the influence function of $T^*_{PRD}$
and the finite sample breakdown point of $T^*_n$  are revealed here as well.
\vs
The competitor $T^*_{RD}$ in RH99 has an advantage over $T^*_{PRD}$ in terms of computation in practice, though both confront a challenging computation problem.  The computing issue of $T^*_{RD}$ has been briefly addressed in RH99 (that of its location counterpart, the halfspace median, has been addressed in Liu, et al (2017), among others).
That of $T^*_{PRD}$ is yet to be thoroughly investigated elsewhere.\vs
$T^*_{PRD}$, on the other hand, is superior to $T^*_{RD}$ in terms of breakdown point robustness and is not inferior to $T^*_{RD}$  in terms of relative efficiency.


\begin{center}
{\textbf{\large Acknowledgments}}
\end{center}
The author thanks Professor Emeritus James Stapleton for his careful English proofreading  and an anonymous referee who provided insightful comments and suggestions 
which have led to significant improvements of the manuscript.

\end{document}